\documentclass[11pt]{amsart}
\usepackage{latexsym}
\usepackage{amsmath}
\usepackage{amssymb}

\newcommand{\eproof}{\mbox{\ }\hfill $\Box$ \par \vskip 10pt}

\newtheorem{Theorem}{Theorem}[section]
\newtheorem{lemma}[Theorem]{Lemma}
\newtheorem{prop}[Theorem]{Proposition}

\newtheorem{corol}[Theorem]{Corollary}

\numberwithin{equation}{section}

\def\cal{\mathcal}

\begin{document}

\title[Interior transmission problems]{Interior transmission problems with coefficients of low regularity}

\author[G. Vodev]{Georgi Vodev}

\address {Universit\'e de Nantes, Laboratoire de Math\'ematiques Jean Leray, 2 rue de la Houssini\`ere, BP 92208, 44322 Nantes Cedex 03, France}
\email{Georgi.Vodev@univ-nantes.fr}

\date{}

\begin{abstract} We obtain parabolic transmission eigenvalue-free regions for both isotropic and anisotropic 
interior transmission problems with $L^\infty$ coefficients which are Lipschitz near the boundary. We also
suppose that the restrictions of the coefficients on the boundary are $C^\mu$ smooth with an integer $\mu\ge 2$ and we investigate the way
in which the transmission eigenvalue-free regions depend on $\mu$. 

\end{abstract} 

\maketitle

\setcounter{section}{0}
\section{Introduction}

Let $\Omega\subset\mathbb{R}^d$, $d\ge 2$, be a bounded, connected domain with a $C^\infty$ smooth boundary $\Gamma=\partial\Omega$, and consider
the interior transmission problem

\begin{equation}\label{eq:1.1}
\left\{
\begin{array}{l}
(\nabla c_1(x)\nabla+\lambda^2n_1(x))u_1=0\quad \mbox{in}\quad\Omega,\\
(\nabla c_2(x)\nabla+\lambda^2n_2(x))u_2=0\quad \mbox{in}\quad\Omega,\\
u_1=u_2,\,c_1\partial_\nu u_1=c_2\partial_\nu u_2 \quad\mbox{on}\quad\Gamma,
\end{array}
\right.
\end{equation}
where $\lambda\in \mathbb{C}$, ${\rm Re}\,\lambda\ge 0$, $\nu$ denotes the Euclidean unit inner normal to $\Gamma$ and $c_j,n_j\in L^\infty(\Omega)$, $j=1,2$, are real-valued functions satisfying $c_j(x)\ge b_0$, $n_j(x)\ge b_0$ for some constant $b_0>0$. 
We also suppose that the coefficients are Lipschitz near the boundary. Given a parameter $0<\delta\ll 1$, set 
$\Omega_\delta=\{x\in\Omega:{\rm dist}(x,\Gamma)<\delta\}$. More precisely, we suppose that 
\begin{equation}\label{eq:1.2}
c_j, n_j\in C^1(\overline\Omega_\delta),\quad j=1,2,
\end{equation}
for some $\delta$. Throughout this paper, given an integer $k\ge 1$, $C^k$ will denote the space of the functions $a$ such that
$\partial_x^\alpha a\in L^\infty$ for all multi-indices $\alpha$ with $|\alpha|=k$. 

In this paper we will consider two types of interior transmission problems. The isotropic one when we have the condition
\begin{equation}\label{eq:1.3}
c_1(x)\equiv c_2(x)\equiv 1\quad\mbox{in}\quad\Omega,\quad n_1(x)\neq n_2(x)\quad\mbox{on}\quad\Gamma,
\end{equation}
and the anisotropic one when we have
\begin{equation}\label{eq:1.4}
(c_1(x)-c_2(x))(c_1(x)n_1(x)-c_2(x)n_2(x))\neq 0\quad\mbox{on}\quad\Gamma.
\end{equation}
If the equation (\ref{eq:1.1})
has a non-trivial solution $(u_1,u_2)$ the complex number $\lambda$ is said to be an interior transmission eigenvalue. 
An important question in the theory of the transmission eigenvalues is that one of finding conditions on the coefficients that guarantee
that the transmission eigenvalues form a discrete set. Various such conditions were found in the papers \cite{kn:BCH}, 
\cite{kn:LV1}, \cite{kn:LV2}, \cite{kn:NN1}, \cite{kn:S1}, the most general ones being those in \cite{kn:NN1}. It follows from the
results in \cite{kn:NN1} that under the above conditions the transmission eigenvalues are discrete in both the isotropic and anisotropic
cases, so it is natural to ask if the counting function of these transmission eigenvalues admits Weyl asymptotics. 
Indeed, such asymptotics were obtained in \cite{kn:PV}, \cite{kn:R} in the case of $C^\infty$ smooth coefficients,
and more recently in \cite{kn:FN}, \cite{kn:NN2} for coefficients of very low regularity but with much worse remainder terms. 
The proof in \cite{kn:PV} relies heavily on the location of the transmission eigenvalues on the complex plane and the $C^\infty$
regularity of the coefficients is not essential. Roughly speaking, the result in \cite{kn:PV} says that parabolic 
eigenvalue-free regions imply Weyl asymptotics with a remainder term depending on the shape of the eigenvalue-free region,
namely, the biger the eigenvalue-free region is, the smaller the remainder term is. In the $C^\infty$ setting, 
parabolic eigenvalue-free regions were obtained in the papers \cite{kn:V1}, \cite{kn:V2}, \cite{kn:V3}
in both the isotropic and anisotropic cases. It follows from the results in \cite{kn:V1} and \cite{kn:V3}
that under the condition (\ref{eq:1.3}) there are no transmission eigenvalues in the region

\begin{equation}\label{eq:1.5}
|{\rm Im}\,\lambda|\ge C
\end{equation}
for some constant $C>0$. Note that the eigenvalue-free region (\ref{eq:1.5}) is optimal as shown in Section 4 of \cite{kn:LC}. 
In the one dimensional case (\ref{eq:1.5}) was obtained in \cite{kn:S2}. In the anisotropic case it is proved 
in \cite{kn:V3} that under the condition 
\begin{equation}\label{eq:1.6}
(c_1(x)-c_2(x))(c_1(x)n_1(x)-c_2(x)n_2(x))<0\quad\mbox{on}\quad\Gamma,
\end{equation}
there are no transmission eigenvalues in the region
\begin{equation}\label{eq:1.7}
{\rm Re}\,\lambda\ge 1,\,\,|{\rm Im}\,\lambda|\ge C.
\end{equation}
On the other hand, it is proved 
in \cite{kn:V1} that under the condition 
\begin{equation}\label{eq:1.8}
(c_1(x)-c_2(x))(c_1(x)n_1(x)-c_2(x)n_2(x))>0\quad\mbox{on}\quad\Gamma,
\end{equation}
there are no transmission eigenvalues in the region
\begin{equation}\label{eq:1.9}
|\lambda|\ge C,\quad |{\rm Im}\,\lambda|\ge C|\lambda|^{3/5}.
\end{equation}
Moreover, if in addition to (\ref{eq:1.8}) the condition
\begin{equation}\label{eq:1.10}
\frac{c_1(x)}{n_1(x)}\neq\frac{c_2(x)}{n_2(x)}\quad\mbox{on}\quad\Gamma,
\end{equation}
is assumed, it is proved in \cite{kn:V3} that there are no transmission eigenvalues in the region 
\begin{equation}\label{eq:1.11}
|\lambda|\ge C,\quad |{\rm Im}\,\lambda|\ge C\log|\lambda|.
\end{equation}

Our goal in the present paper is to extend the above results to coefficients having as low regularity as possible. 
The eigenvalue-free regions we get are not as big as those in the $C^\infty$ case, but they are still parabolic. This fact is important
in view of possible applications to obtaining Weyl asymptotics as in \cite{kn:PV}. Our main result is the following

\begin{Theorem} Suppose that the coefficients satisfy (\ref{eq:1.2}) and that $c_j|_\Gamma, n_j|_\Gamma\in C^\mu(\Gamma)$
with some integer $\mu\ge 2$. 
Then, under the condition (\ref{eq:1.3}), there exists a constant
$C>2$ such that there are no transmission eigenvalues in the region
\begin{equation}\label{eq:1.12}
|\lambda|\ge C,\quad |{\rm Im}\,\lambda|\ge C|\lambda|^{1-p_1}
\left(\log|\lambda|\right)^{p_2}, 
\end{equation}
where
$$p_1=\frac{\mu}{2\mu+2d-1}, \quad p_2=\frac{2}{2\mu+2d-1},\quad\quad if \quad\mu\le 2d-1,$$
$$p_1=\frac{\mu+2}{2\mu+2d+5}, \quad p_2=\frac{2}{2\mu+2d+5},\quad\quad if \quad 2d-1<\mu\le 4d,$$
 and $p_1=\frac{2}{5}$, $p_2=0$ if $\mu>4d$. Under the condition (\ref{eq:1.6}), there are no transmission eigenvalues in the region
 \begin{equation}\label{eq:1.13}
|\lambda|\ge C,\quad |{\rm Im}\,\lambda|\ge C|\lambda|^{1-p_1}
\left(\log|\lambda|\right)^{p_2}, \quad {\rm Re}\,\lambda\ge C|\lambda|^{1-p_3}
\left(\log|\lambda|\right)^{p_4},
\end{equation}
where $p_1$ and $p_2$ are the same as above, and
 $$p_3=\frac{\mu}{2\mu+2d+2}, \quad p_4=\frac{1}{\mu+d+1}.$$
 Under the condition (\ref{eq:1.8}), there are no transmission eigenvalues in the region (\ref{eq:1.12}) with
 $$p_1=\frac{\mu}{2\mu+2d+2}, \quad p_2=\frac{1}{\mu+d+1},\quad\quad if \quad \mu\le \frac{4}{3}(d+1),$$
 and $p_1=\frac{2}{7}$, $p_2=0$ if $\mu>\frac{4}{3}(d+1)$.
\end{Theorem}

\noindent
{\bf Remark.} In fact, the theorem remains valid for all $\mu\ge 2$ not necessairily integers. However,
we prefer assuming that $\mu$ is integer because this simplifies the exposition significantly.

This result seems to be new even in the radial case when $\Omega=\{x\in \mathbb{R}^d:|x|\le 1\}$ and $c_j=c_j(r)$, $n_j=n_j(r)$ depend only on the radial
variable $r=|x|$. Since in this case the restrictions of $c_j$ and $n_j$ on the boundary are constants, the above theorem implies the following

\begin{corol} Suppose that the coefficients $c_j(r)$, $n_j(r)$ are Lipschitz on the interval $[1-\delta,1]$ for some $0<\delta\ll 1$. 
 Then, under the condition (\ref{eq:1.3}), there exists a constant
$C>2$ such that there are no transmission eigenvalues in the region
\begin{equation}\label{eq:1.14}
|\lambda|\ge C,\quad |{\rm Im}\,\lambda|\ge C|\lambda|^{3/5}.
\end{equation}
Under the condition (\ref{eq:1.6}), there are no transmission eigenvalues in the region
 \begin{equation}\label{eq:1.15}
|\lambda|\ge C,\quad |{\rm Im}\,\lambda|\ge C|\lambda|^{3/5}, \quad {\rm Re}\,\lambda\ge C_\epsilon|\lambda|^{1/2+\epsilon},
\end{equation}
for every $0<\epsilon\ll 1$. 
 Under the condition (\ref{eq:1.8}), there are no transmission eigenvalues in the region 
 \begin{equation}\label{eq:1.16}
|\lambda|\ge C,\quad |{\rm Im}\,\lambda|\ge C|\lambda|^{5/7}.
\end{equation}
\end{corol}

To prove Theorem 1.1 we adapt the methods developed in \cite{kn:V1}, \cite{kn:V2}, \cite{kn:V3} for $C^\infty$ smooth coefficients. 
In these papers the transmission eigenvalue-free regions were derived from suitable approximations of the interior Dirichlet-to-Neumann map
outside a parabolic neighbourhood of the real axis. In other words, the problem of finding transmission eigenvalue-free regions
was reduced to that one of finding as good as possible approximation of the interior Dirichlet-to-Neumann map by a semiclassical 
pseudodifferential operator ($h-\Psi$DO) on the boundary (with a semiclassical parameter $h\sim |\lambda|^{-1}$) outside as small as possible parabolic neighbourhood of the real axis.
With such an approximation in hands, in order to get the eigenvalue-free regions, one has to invert a semiclassical 
pseudodifferential operator whose symbol, say $a$, can be camputed explicitly in terms of the restrictions of the coefficients
on the boundary and the principal symbol of the Laplace-Beltrami operator on the boundary. Recall that the boundary $\Gamma$
can be considered as a compact Riemannian manifold of dimension $d-1$ without boundary with a Riemannian metric
induced by the Euclidean one. Then the conditions (\ref{eq:1.3}) and (\ref{eq:1.4}) guarantee that the function $a$
is invertible outside some parabolic neighbourhood of the real axis. Then to show that the semiclassical 
pseudodifferential operator ${\rm Op}_h(a)$ with symbol $a$ is invertible, one has to make use of some pseudodifferential calculas which are  well-known in the $C^\infty$ setting. Thus one arrives at the conclusion that the eigenvalue-free regions correspond to the regions
where the Dirichlet-to-Neumann map has a good approximation and where the operator ${\rm Op}_h(a)$ is invertible. 
Therefore, to get as big as possible eigenvalue-free regions one has to find as good as possible approximation of the
 interior Dirichlet-to-Neumann map. In the $C^\infty$ case this is done in \cite{kn:V1}, \cite{kn:V2}, \cite{kn:V3}
 by constructing semiclassical parametrices for the solutions of the equation (\ref{eq:1.1}) 
 near the boundary mod ${\cal O}(h^\infty)$. Note that
 the $C^\infty$ regularity of the coefficients near the boundary is essential in the analysis in these papers. 
 The main difficulty in constructing such parametrices is that one has to solve the eikonal equation mod ${\cal O}(x_1^\infty)$
 and the transport equations
 mod ${\cal O}(x_1^\infty+h^\infty)$, where $0<x_1\ll 1$ is the normal coordinate near the boundary, 
 that is the distance to $\Gamma$. 
 This is especially difficult to do when the boundary data is microlocally supported near the glancing region.
 Of course, when the coefficients are of low regularity this analysis does not work any more. In this case we built (see Section 4)
 a less accurate parametrix by solving the eikonal equation mod ${\cal O}(x_1)$ and with no need to solve the transport equations.
 In other words, our parametrix is the simplest possible in this context. Nevertheless, this leads to some approximation of 
 the Dirichlet-to-Neumann map, which of course is not as good as in the $C^\infty$ case. Consequently, the 
 eigenvalue-free regions we get by using this approximation are much smaller than those in the $C^\infty$ case.
 In order to make this approach work, however, we need to use $h-\Psi$DOs with symbols which are $C^\mu$ smooth 
 with respect to the space variable. To this end, we develop in Section 2 some pseudodifferential calculas
 for such operators. In particular, we find some useful criteria for $L^2$ boundedness 
 (see Proposition 2.4) and we also show how to invert such operators (see Proposition 2.6).
 Note finally that $\mu=2$ is the lowest regularity that makes possible the construction of some parametrix. 

The paper is organized as follows. In Section 2 we recall some basic properties of the $h-\Psi$DOs with $C^\infty$ symbols
concerning the $L^2$ boundedness and the composition of two operators. We then extend these properties to 
$h-\Psi$DOs with symbols which are $C^\mu$ smooth with respect to the space variable. This is done by using
a suitable interpolation between symbols which are $L^\infty$ and $C^\infty$ smooth with respect to the space variable.
To this end, we establish a criteria for $L^2$ boundedness for $h-\Psi$DOs with $L^\infty$ 
symbols (see Proposition 2.3). In Section 3 we prove a priori estimates
which allow us to bound the norm of the difference between the Dirichlet-to-Neumann map and the parametrix. 
In Section 4 we build a parametrix for the solutions of the equation (\ref{eq:1.1}) near the boundary, which is a finite sum of $h-$FIOs.  
We then use it to get the parametrix for the Dirichlet-to-Neumann map.
In Section 5 we improve our parametrix in the elliptic region in the isotropic case. This allows us to get a better approximation
of the Dirichlet-to-Neumann map in this case. This improvement is crussial in order to get the eigenvalue-free regions
in the isotropic case. In Section 6 we invert 
the operator ${\rm Op}_h(a)$ with a symbol $a$ as above outside some parabolic neighbourhoods of the real axis
as well as outside some parabolic neighbourhoods of the imaginary axis when the condition (\ref{eq:1.6}) is assumed.
This gives us the eigenvalue-free regions. 
Note that in our case $a$ 
 is $C^\mu$ smooth with respect to the space variable, so we need to use the pseudodifferential calculas developed in Section 2.

 \section{$h-\Psi$DOs with symbols of low regularity} 
 
 We begin this section by recalling some basic properties of the $h-\Psi$DOs with $C^\infty$ symbols. In what follows we will
 identify the cotangent space $T^*\mathbb{R}^{d-1}$ with $\mathbb{R}^{d-1}\times \mathbb{R}^{d-1}$.  Introduce
 the space $S^k$, $k\in \mathbb{R}$, of all functions $a(x,\xi)\in C^\infty(\mathbb{R}^{d-1}\times \mathbb{R}^{d-1})$ satisfying
\begin{equation}\label{eq:2.1}
\left|\partial_\xi^\alpha \partial_x^\beta a(x,\xi)\right|\le C_{\alpha,\beta}\langle\xi\rangle^{k-|\alpha|}
\end{equation}
for all multi-indices $\alpha$ and $\beta$, where $\langle\xi\rangle:=(1+|\xi|^2)^{1/2}$.
We define 
the $h-\Psi$DO with a symbol $a\in S^k$ by
$$\left({\rm Op}_h(a)f\right)(x)=(2\pi h)^{-d+1}\int_{\mathbb{R}^{d-1}}\int_{\mathbb{R}^{d-1}} e^{-\frac{i}{h}\langle x-y,\xi\rangle}a(x,\xi)f(y)d\xi dy,$$
where $0<h\le 1$ is a semiclassical parameter. 
Denote by $H_h^k(\mathbb{R}^{d-1})$ the Sobolev space equipped with the $h$-semiclassical norm
$$\left\|u\right\|_{H_h^k}:= \left\|{\rm Op}_h(\langle\xi\rangle^k)u\right\|_{L^2},$$ 
where $\|\cdot\|_{L^2}$ denotes the norm in $L^2(\mathbb{R}^{d-1})$. 
It is well-known that the operator ${\rm Op}_h(a):H_h^k(\mathbb{R}^{d-1})\to L^2(\mathbb{R}^{d-1})$ is bounded
for symbols $a\in S^k$. In fact, we have a stronger result.

\begin{prop} Let $a\in C^\infty(\mathbb{R}^{d-1}\times \mathbb{R}^{d-1})$ satisfy 
\begin{equation}\label{eq:2.2}
\left|\partial_\xi^\alpha \partial_x^\beta a(x,\xi)\right|\le C_{\alpha,\beta}\langle\xi\rangle^k
\end{equation}
for all multi-indices $\alpha$ and $\beta$.
 Then there exists an integer 
$s_d$ depending only on the dimension such that we have the bound
\begin{equation}\label{eq:2.3}
\left\|{\rm Op}_h(a)\right\|_{H_h^k\to L^2}
\lesssim \sum_{|\alpha|+|\beta|\le s_d}h^{\frac{|\alpha|+|\beta|}{2}}C_{\alpha,\beta}.
\end{equation}
\end{prop}

{\it Proof.} Since ${\rm Op}_h(a)={\rm Op}_h(a\langle\xi\rangle^{-k}){\rm Op}_h(\langle\xi\rangle^k)$ and the function
$a\langle\xi\rangle^{-k}$ satisfies (\ref{eq:2.2}) with $k=0$, it suffices to prove (\ref{eq:2.3}) for $k=0$. In this case, 
it is easy to see that the norm in the left-hand side of (\ref{eq:2.3}) is equal to the norm
$\left\|{\rm Op}_1(a_h)\right\|_{L^2\to L^2}$, where $a_h(x,\xi)=a(h^{1/2}x,h^{1/2}\xi)$. 
On the other hand, by Theorem 18.6.3 of \cite{kn:H} (see also Theorem 7.11 of \cite{kn:DS}) the operator ${\rm Op}_1(a_h):L^2\to L^2$ is bounded and
$$\left\|{\rm Op}_1(a_h)\right\|_{L^2\to L^2}\lesssim\sum_{|\alpha|+|\beta|\le s_d}
{\rm sup}\left|\partial_\xi^\alpha \partial_x^\beta a_h(x,\xi)\right|,$$
which clearly implies (\ref{eq:2.3}) with $k=0$.
\eproof

We will now use this proposition to prove the following

\begin{prop} Let $a_1$ satisfy (\ref{eq:2.2}) with $k=0$ and constants $C_{\alpha,\beta}^{(1)}$, and 
let $a_2$ satisfy (\ref{eq:2.2}) with constants $C_{\alpha,\beta}^{(2)}$. 
 Then we have the bound
\begin{equation}\label{eq:2.4}
\left\|{\rm Op}_h(a_1){\rm Op}_h(a_2)-{\rm Op}_h(a_1a_2)\right\|_{H_h^k\to L^2}$$
 $$\lesssim \sum_{|\alpha_1|+|\beta_1|+|\alpha_2|+|\beta_2|\le s'_d,\,|\alpha_1|+|\beta_1|\ge 1,\,|\alpha_2|+|\beta_2|\ge 1}h^{\frac{|\alpha_1|+|\beta_1|+|\alpha_2|+|\beta_2|}{2}}C^{(1)}_{\alpha_1,\beta_1}C^{(2)}_{\alpha_2,\beta_2},
\end{equation}
where $s'_d=2d+5+s_d$.
\end{prop}

{\it Proof.} We will use the fact that 
$${\rm Op}_h(a_1){\rm Op}_h(a_2)={\rm Op}_h(b),$$
where
$$b(x,\xi)=e^{ih\langle D_\eta,D_y\rangle}a_1(x,\eta)a_2(y,\xi)|_{y=x,\eta=\xi}$$
$$=(2\pi)^{-d+1}\int_{\mathbb{R}^{d-1}}\int_{\mathbb{R}^{d-1}} e^{i\langle\zeta_1,\zeta_2\rangle}a_1(x,\xi-h^{1/2}\zeta_1)a_2(x-h^{1/2}\zeta_2,\xi)d\zeta_1 d\zeta_2,$$
where we have put $D_\eta=-i\nabla_\eta$, $D_y=-i\nabla_y$. 
Set $\zeta=(\zeta_1,\zeta_2)\in \mathbb{R}^{2d-2}$ and $\varphi(\zeta)=\langle\zeta_1,\zeta_2\rangle$. Clearly,
$$\nabla_\zeta\varphi=(\nabla_{\zeta_1}\varphi,\nabla_{\zeta_2}\varphi)=(\zeta_2,\zeta_1),$$
so we have $|\nabla_\zeta\varphi|=|\zeta|$. 
Consider now the function
$$g(t)=a_1(x,\xi-th^{1/2}\zeta_1)a_2(x-th^{1/2}\zeta_2,\xi),\quad 0\le t\le 1.$$
An easy computation leads to the formulas
$$g'(t)=-h^{1/2}\langle\zeta_1,\nabla_\xi a_1(x,\xi-th^{1/2}\zeta_1)\rangle a_2(x-th^{1/2}\zeta_2,\xi)$$ 
$$-h^{1/2}\langle\zeta_2,\nabla_x a_2(x-th^{1/2}\zeta_2,\xi)\rangle a_1(x,\xi-th^{1/2}\zeta_1),$$
$$g''(t)={\cal G}_1(t)+{\cal G}_2(t),$$
where
$${\cal G}_1=h\langle\zeta_1,\nabla_\xi\langle\zeta_1,\nabla_\xi a_1(x,\xi-th^{1/2}\zeta_1)\rangle\rangle a_2(x,\xi)$$ 
$$+h\langle\zeta_2,\nabla_x\langle\zeta_2,\nabla_x a_2(x-th^{1/2}\zeta_2,\xi)\rangle\rangle a_1(x,\xi),$$
$${\cal G}_2=h\langle\zeta_1,\nabla_\xi\langle\zeta_1,\nabla_\xi a_1(x,\xi-th^{1/2}\zeta_1)\rangle\rangle
\left(a_2(x-th^{1/2}\zeta_2,\xi)-a_2(x,\xi)\right)$$ 
$$+h\langle\zeta_2,\nabla_x\langle\zeta_2,\nabla_x a_2(x-th^{1/2}\zeta_2,\xi)\rangle\rangle
\left(a_1(x,\xi-th^{1/2}\zeta_1)-a_1(x,\xi)\right)$$
$$+2h\langle\zeta_1,\nabla_\xi a_1(x,\xi-th^{1/2}\zeta_1)\rangle\langle\zeta_2,\nabla_x a_2(x-th^{1/2}\zeta_2,\xi)\rangle.$$
Clearly, the function ${\cal G}_2$ satisfies the bounds
\begin{equation}\label{eq:2.5}
\left|\partial_\zeta^\gamma{\cal G}_2\right|\lesssim\langle\xi\rangle^k\sum_{\ell=0}^3|\zeta|^{\ell}
\sum_{|\alpha|+|\beta|=|\gamma|+\ell,\,|\alpha|\ge 1,\,|\beta|\ge 1}
h^{\frac{|\alpha|+|\beta|}{2}}C^{(1)}_{\alpha,0}C^{(2)}_{0,\beta}
\end{equation}
for all multi-indices $\gamma$, uniformly in $t$. Observe now that we can write 
$$g(1)=g(0)+g'(0)+\int_0^1(1-t)g''(t)dt,$$
where
$$g(0)=a_1(x,\xi)a_2(x,\xi),$$
$$g'(0)=-h^{1/2}\langle\zeta_1,\nabla_\xi a_1(x,\xi)\rangle a_2(x,\xi)
-h^{1/2}\langle\zeta_2,\nabla_x a_2(x,\xi)\rangle a_1(x,\xi).$$
On the other hand, given any function $\psi\in C^\infty(\mathbb{R}^{d-1})$, we have the formula
$$(2\pi)^{-d+1}\int_{\mathbb{R}^{2d-2}}e^{i\varphi(\zeta)}\psi(\zeta_j)d\zeta=\psi(0)$$
where $j=1,2$. Therefore we have
$$\int_{\mathbb{R}^{2d-2}}e^{i\varphi}g'(0)d\zeta=\int_{\mathbb{R}^{2d-2}}e^{i\varphi}{\cal G}_1(t)d\zeta=0.$$
We get from the above formulas
\begin{equation}\label{eq:2.6}
b-a_1a_2=(2\pi)^{-d+1}\int_0^1(1-t)\int_{\mathbb{R}^{2d-2}}e^{i\varphi}{\cal G}_2(t)d\zeta dt.
\end{equation}
Let $\phi\in C_0^\infty(\mathbb{R})$,
$\phi(\sigma)=1$ for $|\sigma|\le 1$, $\phi(\sigma)=0$ for $|\sigma|\ge 2$. We will now bound the integrals
$$I_1=(2\pi)^{-d+1}\int_{\mathbb{R}^{2d-2}}e^{i\varphi}\phi(|\zeta|){\cal G}_2(t)d\zeta,$$
$$I_2=(2\pi)^{-d+1}\int_{\mathbb{R}^{2d-2}}e^{i\varphi}(1-\phi)(|\zeta|){\cal G}_2(t)d\zeta.$$
By (\ref{eq:2.5}) with $\gamma=0$, we have
\begin{equation}\label{eq:2.7}
\left|I_1\right|\lesssim\langle\xi\rangle^k\sum_{|\alpha|+|\beta|\le 3,\,|\alpha|\ge 1,\,|\beta|\ge 1}
h^{\frac{|\alpha|+|\beta|}{2}}C^{(1)}_{\alpha,0}C^{(2)}_{0,\beta}
\end{equation}
 uniformly in $t$. To bound the second integral we will integrate by parts. To this end observe that
 $$L_\zeta e^{i\varphi}=e^{i\varphi},$$
 where 
 $$L_\zeta =-i|\zeta|^{-2}\langle\nabla_\zeta\varphi,\nabla_\zeta\rangle.$$
 It is easy to see that, given any integer $m\ge 1$, the adjoint operator to $L_\zeta^m$ satisfies
 \begin{equation}\label{eq:2.8}
\left|(L_\zeta^m)^*f(\zeta)\right|\lesssim |\zeta|^{-m}\sum_{|\gamma|\le m}\left|\partial_\zeta^\gamma f(\zeta)\right|.
\end{equation}
We can now write the integral $I_2$ in the form
$$I_2=(2\pi)^{-d+1}\int_{\mathbb{R}^{2d-2}}e^{i\varphi}(L_\zeta^m)^*\left((1-\phi){\cal G}_2(t)\right)d\zeta.$$
Using (\ref{eq:2.5}) and (\ref{eq:2.8}) we obtain
 \begin{equation}\label{eq:2.9}
|I_2|\lesssim\langle\xi\rangle^k\sum_{\ell=0}^3\int_{\zeta\in\mathbb{R}^{2d-2}:|\zeta|\ge 1}|\zeta|^{\ell-m}d\zeta
\sum_{|\alpha|+|\beta|\le m+\ell,\,|\alpha|\ge 1,\,|\beta|\ge 1}
h^{\frac{|\alpha|+|\beta|}{2}}C^{(1)}_{\alpha,0}C^{(2)}_{0,\beta}$$
$$\lesssim\langle\xi\rangle^k
\sum_{|\alpha|+|\beta|\le 2d+5,\,|\alpha|\ge 1,\,|\beta|\ge 1}
h^{\frac{|\alpha|+|\beta|}{2}}C^{(1)}_{\alpha,0}C^{(2)}_{0,\beta}
\end{equation}
if we take $m=2d+2$. By (\ref{eq:2.6}), (\ref{eq:2.7}) and (\ref{eq:2.9}), we conclude
\begin{equation}\label{eq:2.10}
|b-a_1a_2|\lesssim\langle\xi\rangle^k
\sum_{|\alpha_1|+|\beta_2|\le 2d+5,\,|\alpha_1|\ge 1,\,|\beta_2|\ge 1}
h^{\frac{|\alpha_1|+|\beta_2|}{2}}C^{(1)}_{\alpha_1,0}C^{(2)}_{0,\beta_2}.
\end{equation}
Similarly, given any multi-indices $\alpha$ and $\beta$, we can obtain the bounds
\begin{equation}\label{eq:2.11}
\langle\xi\rangle^{-k}h^{\frac{|\alpha|+|\beta|}{2}}|\partial_\xi^\alpha\partial_x^\beta(b-a_1a_2)|$$
$$\lesssim
\sum_{|\alpha_1|+|\alpha_2|+|\beta_1|+|\beta_2|\le |\alpha|+|\beta|+2d+5,\,|\alpha_1|+|\beta_1|\ge 1,\,|\alpha_2|+|\beta_2|\ge 1}
h^{\frac{|\alpha_1|+|\alpha_2|+|\beta_1|+|\beta_2|}{2}}C^{(1)}_{\alpha_1,\beta_1}C^{(2)}_{\alpha_2,\beta_2}.
\end{equation}
Clearly, the bound (\ref{eq:2.4}) follows from (\ref{eq:2.11}) and Proposition 2.1.
\eproof

The above definition of the $h-\Psi$DOs makes sense also for symbols $a(x,\xi)$ which are smooth
only with respect to the variable $\xi$ and satisfy
\begin{equation}\label{eq:2.12}
\left|\partial_\xi^\alpha a(x,\xi)\right|\le C_\alpha\langle\xi\rangle^{k-|\alpha|}
\end{equation}
for all multi-indices $\alpha$, where $k\in \mathbb{R}$. For such operators, we will prove the following

\begin{prop} Let $a$ satisfy (\ref{eq:2.12}) and let $k'>k$. 
  Then we have
\begin{equation}\label{eq:2.13}
\left\|{\rm Op}_h(a)\right\|_{H_h^{k'}\to L^2}\lesssim \log(1+h^{-1})\sum_{|\alpha|\le d}C_\alpha.
\end{equation}
\end{prop}

{\it Proof.} As above, it suffices to prove (\ref{eq:2.13}) for $k'=0$. Then we have $k<0$. Clearly, the kernel $K$ of the operator 
${\rm Op}_h(a)$ is
$$K(x,y)=(2\pi h)^{-d+1}\int_{\mathbb{R}^{d-1}} e^{-\frac{i}{h}\langle x-y,\xi\rangle}a(x,\xi)d\xi.$$
In view of Schur's lemma, it suffices to bound from above the integrals
$$\int_{\mathbb{R}^{d-1}}|K(x,y)|dx+\int_{\mathbb{R}^{d-1}}|K(x,y)|dy.$$
Integrating by parts in the above integral leads to the bounds
$$|K(x,y)|\lesssim h^{-d+1}\left(\frac{h}{|x-y|}\right)^{d}\sum_{|\alpha|=d}
\int_{\mathbb{R}^{d-1}}\left|\partial_\xi^\alpha a(x,\xi)\right|d\xi$$
 $$\lesssim\frac{h}{|x-y|^{d}}\sum_{|\alpha|=d}C_\alpha\int_{\mathbb{R}^{d-1}}\langle\xi\rangle^{k-d}d\xi$$
 $$\lesssim\frac{h}{|x-y|^{d}}\sum_{|\alpha|=d}C_\alpha.$$
 Hence
 $$\int_{|x-y|\ge h}|K(x,y)|dx+\int_{|x-y|\ge h}|K(x,y)|dy$$ 
 $$\lesssim h\sum_{|\alpha|=d}C_\alpha\int_{w\in \mathbb{R}^{d-1}:|w|\ge h}|w|^{-d}dw$$
$$\lesssim h\sum_{|\alpha|=d}C_\alpha\int_h^\infty\sigma^{-2}d\sigma\lesssim \sum_{|\alpha|=d}C_\alpha.$$
We will now bound the kernel for $|x-y|<h$. Set $A=\frac{h}{|x-y|}>1$. Let the function $\phi$ be as above and 
decompose the kernel as $K=K_1+K_2$, where
$$K_1(x,y)=(2\pi h)^{-d+1}\int_{\mathbb{R}^{d-1}} e^{-\frac{i}{h}\langle x-y,\xi\rangle}a(x,\xi)\phi(|\xi|/A)d\xi.$$
An integration by parts leads to the bounds
$$|K_1(x,y)|\lesssim h^{-d+1}\left(\frac{h}{|x-y|}\right)^{d-2}\sum_{|\alpha|=d-2}
\int_{\mathbb{R}^{d-1}}\left|\partial_\xi^\alpha\left(a(x,\xi)\phi(|\xi|/A)\right)\right|d\xi$$
$$\lesssim \frac{h^{-1}}{|x-y|^{d-2}}\sum_{|\alpha|\le d-2}A^{|\alpha|+2-d}
\int_{|\xi|\le A}\left|\partial_\xi^\alpha a(x,\xi)\right|d\xi$$
$$\lesssim \frac{h^{-1}}{|x-y|^{d-2}}\sum_{|\alpha|\le d-2}
C_\alpha A^{|\alpha|+2-d}\int_{|\xi|\le A}\langle\xi\rangle^{k-|\alpha|}d\xi$$
$$\lesssim \frac{h^{-1}}{|x-y|^{d-2}}\sum_{|\alpha|\le d-2}
C_\alpha A^{|\alpha|+2-d}\int_0^{A+1}\sigma^{k-|\alpha|+d-2}d\sigma$$
$$\lesssim \frac{h^{-1}A^{1+k}}{|x-y|^{d-2}}\sum_{|\alpha|\le d-2}C_\alpha$$
$$\lesssim \frac{A^{k}}{|x-y|^{d-1}}\sum_{|\alpha|\le d-2}C_\alpha.$$
Set $\epsilon=\left(\log(1+h^{-1})\right)^{-1}$. Clearly, there is $0<h_0\ll 1$ such that $\epsilon<-k$ for $0<h<h_0$, and hence $A^k<A^{-\epsilon}$.
Thus, taking into account that $h^{-\epsilon}=e$, we get
$$|K_1(x,y)|\lesssim\frac{1}{|x-y|^{d-1-\epsilon}}\sum_{|\alpha|\le d-2}C_\alpha$$
for $0<h<h_0$. Clearly, for $h_0\le h\le 1$ the above bound still holds with $0<\epsilon\ll 1$ independent of $h$. 
Hence, in all cases, we have 
 $$\int_{|x-y|\le h}|K_1(x,y)|dx+\int_{|x-y|\le h}|K_1(x,y)|dy$$ 
 $$\lesssim \sum_{|\alpha|\le d-2}C_\alpha\int_{w\in\mathbb{R}^{d-1}:|w|\le h}|w|^{-d+1+\epsilon}dw$$
$$\lesssim \sum_{|\alpha|\le d-2}C_\alpha\int_0^h\sigma^{-1+\epsilon}d\sigma\lesssim \epsilon^{-1}\sum_{|\alpha|\le d-2}C_\alpha.$$
Similarly, we have 
$$|K_2(x,y)|\lesssim h^{-d+1}\left(\frac{h}{|x-y|}\right)^{d-1}\sum_{|\alpha|=d-1}
\int_{\mathbb{R}^{d-1}}\left|\partial_\xi^\alpha\left(a(x,\xi)(1-\phi)(|\xi|/A)\right)\right|d\xi$$
$$\lesssim \frac{1}{|x-y|^{d-1}}\sum_{|\alpha|\le d-1}A^{|\alpha|+1-d}
\int_{|\xi|\ge A}\left|\partial_\xi^\alpha a(x,\xi)\right|d\xi$$
$$\lesssim \frac{1}{|x-y|^{d-1}}\sum_{|\alpha|\le d-1}
C_\alpha A^{|\alpha|+1-d}\int_{|\xi|\ge A}\langle\xi\rangle^{k-|\alpha|}d\xi$$
$$\lesssim \frac{A^{-\epsilon}}{|x-y|^{d-1}}\sum_{|\alpha|\le d-1}
C_\alpha \int_{\mathbb{R}^{d-1}}\langle\xi\rangle^{k-d+1+\epsilon}d\xi$$
$$\lesssim \frac{1}{|x-y|^{d-1-\epsilon}}\sum_{|\alpha|\le d-1}
C_\alpha \int_0^\infty(\sigma+1)^{k-1+\epsilon}d\sigma$$
$$\lesssim \frac{1}{|x-y|^{d-1-\epsilon}}\sum_{|\alpha|\le d-1}C_\alpha.$$
As above, this implies 
$$\int_{|x-y|\le h}|K_2(x,y)|dx+\int_{|x-y|\le h}|K_2(x,y)|dy\lesssim \epsilon^{-1}\sum_{|\alpha|\le d-1}C_\alpha,$$ 
which clearly completes the proof of (\ref{eq:2.13}) with $k'=0$. 
\eproof
 
 The logarithmic term in the right-hand side of (\ref{eq:2.13}) can be removed if some regularity of the symbol in $x$ is assumed. 
 Indeed, in this case we can decompose $a$ as a sum of a large symbol which
 is $C^\infty$ in $x$ and a small symbol which is $L^\infty$ in $x$. Then we can apply (\ref{eq:2.3}) to the $C^\infty$ part
 and (\ref{eq:2.13}) to the $L^\infty$ part. More precisely,
 we have the following
 
 \begin{prop} Suppose that $a$ satisfy (\ref{eq:2.1}) for all multi-indices
 $\alpha$ and $\beta$ such that $|\beta|\le \mu$, $\mu\ge 1$ being integer, and let $k'>k$. 
  Then we have
\begin{equation}\label{eq:2.14}
\left\|{\rm Op}_h(a)\right\|_{H_h^{k'}\to L^2}\lesssim \sum_{|\alpha|+|\beta|\le s_d,\,|\beta|\le\mu}h^{\frac{|\alpha|+|\beta|}{2}}C_{\alpha,\beta}$$ 
$$+h^{\mu/2}\log(1+h^{-1})\sum_{|\alpha|\le d,\,|\beta|\le\mu}C_{\alpha,\beta}.
\end{equation}
Moreover, if $\mu\ge s_d$ the second sum in the right-hand side of (\ref{eq:2.14}) can be removed. 
\end{prop}

{\it Proof.} Let $\phi_0\in C_0^\infty(\mathbb{R}^{d-1})$, $\phi_0\ge 0$, be such that $\int_{\mathbb{R}^{d-1}}\phi_0(x)dx=1$.
Given a parameter $0<t\le 1$, set
$$a_t(x,\xi)=t^{-d+1}\int_{\mathbb{R}^{d-1}}\phi_0((y-x)/t)a(y,\xi)dy=
\int_{\mathbb{R}^{d-1}}\phi_0(w)a(x+tw,\xi)dw.$$
Clearly, the function $a_t$ is $C^\infty$ smooth in $x$ and satisfies (\ref{eq:2.1}) uniformly in $t$ for all multi-indices
 $\alpha$ and $\beta$ such that $|\beta|\le \mu$, while for $|\beta|\ge \mu+1$ it satisfies the bounds
\begin{equation}\label{eq:2.15}
\left|\partial_\xi^\alpha \partial_x^\beta a_t(x,\xi)\right|\le t^{-|\beta|+\mu}\langle\xi\rangle^{k-|\alpha|}
\sum_{|\beta'|=\mu}C_{\alpha,\beta'}.
\end{equation}
The Taylor expansion of the function $a(x+tw,\xi)$ at $t=0$ yields
$$a(x+tw,\xi)=\sum_{s=0}^{\mu-1}\frac{t^s}{s!}\partial_t^sa(x+tw,\xi)|_{t=0}+\frac{t^\mu}{\mu!}\partial_t^\mu a(x+tw,\xi)|_{t=t'}$$
$$=\sum_{s=0}^{\mu-1}\frac{t^s}{s!}\sum_{|\beta|=s}w^\beta\partial_x^\beta a(x,\xi)+\frac{t^\mu}{\mu!}
 \sum_{|\beta|=\mu}w^\beta\partial_x^\beta a(x+t'w,\xi)$$
with some $0\le t'\le t$. Clearly, the functions $\partial_x^\beta a(x+t'w,\xi)$ with $|\beta|=\mu$ satisfy (\ref{eq:2.12}) uniformly
in $t'w$. Therefore, by Proposition 2.3 we have
\begin{equation}\label{eq:2.16}
\sum_{|\beta|=\mu}\left\|{\rm Op}_h(\partial_x^\beta a(x+t'w,\xi))\right\|_{H_h^{k'}\to L^2}\lesssim \log(1+h^{-1})\sum_{|\alpha|\le d,|\beta|=\mu}C_{\alpha,\beta}.
\end{equation}
We now apply Proposition 2.1 with $a_t$. In view of (\ref{eq:2.15}) we get
\begin{equation}\label{eq:2.17}
\left\|{\rm Op}_h(a_t)\right\|_{H_h^k\to L^2}
\lesssim \sum_{|\alpha|+|\beta|\le s_d,\,|\beta|\le\mu}h^{\frac{|\alpha|+|\beta|}{2}}C_{\alpha,\beta}
\end{equation}
if we take $t=h^{1/2}$. Thus, by (\ref{eq:2.16}) and (\ref{eq:2.17}), we obtain
\begin{equation}\label{eq:2.18}
\left\|{\rm Op}_h(a)\right\|_{H_h^{k'}\to L^2}
\lesssim \sum_{s=1}^{\mu-1}h^{s/2}\sum_{|\beta|=s}\left\|{\rm Op}_h(\partial_x^\beta a)\right\|_{H_h^{k'}\to L^2}$$ 
$$+h^{\mu/2}\log(1+h^{-1})\sum_{|\alpha|\le d,|\beta|=\mu}C_{\alpha,\beta}
+\sum_{|\alpha|+|\beta|\le s_d,\,|\beta|\le\mu}h^{\frac{|\alpha|+|\beta|}{2}}C_{\alpha,\beta},
\end{equation}
where the first sum in the righ-hand side is zero if $\mu=1$. It is easy to see now that (\ref{eq:2.14}) follows from
(\ref{eq:2.18}) by induction in $\mu$. 
\eproof

We will next extend Proposition 2.2 to $h-\Psi$DOs with $C^\mu$ smooth symbols, where $\mu\ge 1$ is integer. 
We will first prove the following

\begin{prop} Let $a_1$ satisfy (\ref{eq:2.1}) with $k<0$ for all multi-indices
$\alpha$ and $\beta$ such that $|\beta|\le\mu$ with constants $C_{\alpha,\beta}^{(1)}$, and 
let $a_2$ satisfy (\ref{eq:2.1}) for all multi-indices
$\alpha$ and $\beta$ with constants $C_{\alpha,\beta}^{(2)}$. 
 Then we have the bound
\begin{equation}\label{eq:2.19}
\left\|{\rm Op}_h(a_1){\rm Op}_h(a_2)-{\rm Op}_h(a_1a_2)\right\|_{H_h^k\to L^2}$$
 $$\lesssim \sum_{|\alpha_1|+|\beta_1|+|\alpha_2|+|\beta_2|\le s'_d,\,|\beta_1|\le\mu,\,|\alpha_1|+|\beta_1|\ge 1,\,|\alpha_2|+|\beta_2|\ge 1}h^{\frac{|\alpha_1|+|\beta_1|+|\alpha_2|+|\beta_2|}{2}}C^{(1)}_{\alpha_1,\beta_1}C^{(2)}_{\alpha_2,\beta_2}$$
 $$+h^{\mu/2}\log(1+h^{-1})\sum_{|\alpha_1|\le d,\,|\beta_1|=\mu,\,|\alpha_2|+|\beta_2|\le s_d}h^{\frac{|\alpha_2|+|\beta_2|}{2}}C^{(1)}_{\alpha_1,\beta_1}C^{(2)}_{\alpha_2,\beta_2}$$ 
 $$+h^{\mu/2}\log(1+h^{-1})\sum_{|\alpha_1|+|\alpha_2|\le d,\,|\beta_1|=\mu,\,\beta_2=0}C^{(1)}_{\alpha_1,\beta_1}C^{(2)}_{\alpha_2,\beta_2}.
\end{equation}
Let $a_1$ satisfy (\ref{eq:2.1}) with $k=0$ for all multi-indices
$\alpha$ and $\beta$ with constants $C_{\alpha,\beta}^{(1)}$, and 
let $a_2$ satisfy (\ref{eq:2.1}) with $k<k'$ for all multi-indices
$\alpha$ and $\beta$ such that $|\beta|\le\mu$ with constants $C_{\alpha,\beta}^{(2)}$. 
 Then we have the bound
\begin{equation}\label{eq:2.20}
\left\|{\rm Op}_h(a_1){\rm Op}_h(a_2)-{\rm Op}_h(a_1a_2)\right\|_{H_h^{k'}\to L^2}$$
 $$\lesssim \sum_{|\alpha_1|+|\beta_1|+|\alpha_2|+|\beta_2|\le s'_d,\,|\beta_2|\le\mu,\,|\alpha_1|+|\beta_1|\ge 1,\,|\alpha_2|+|\beta_2|\ge 1}h^{\frac{|\alpha_1|+|\beta_1|+|\alpha_2|+|\beta_2|}{2}}C^{(1)}_{\alpha_1,\beta_1}C^{(2)}_{\alpha_2,\beta_2}$$
 $$+h^{\mu/2}\log(1+h^{-1})\sum_{|\alpha_1|+|\beta_1|\le s_d,\,|\alpha_2|\le d,\,|\beta_2|=\mu}h^{\frac{|\alpha_1|+|\beta_1|}{2}}C^{(1)}_{\alpha_1,\beta_1}C^{(2)}_{\alpha_2,\beta_2}$$
 $$+h^{\mu/2}\log(1+h^{-1})\sum_{|\alpha_1|+|\alpha_2|\le d,\,\beta_1=0,\,|\beta_2|=\mu}C^{(1)}_{\alpha_1,\beta_1}C^{(2)}_{\alpha_2,\beta_2}.
\end{equation}
Let $a_1$ satisfy (\ref{eq:2.1}) with $k=0$ for all multi-indices
$\alpha$ and $\beta$ with constants $C_{\alpha,\beta}^{(1)}$, and 
let $a_2$ be independent of the variable $\xi$ and satisfies
\begin{equation}\label{eq:2.21}
\left|\partial_x^\beta a_2(x)\right|\le C^{(2)}_{\beta}
\end{equation}
for all multi-indices $\beta$ such that $|\beta|\le\mu$.
 Then we have the bound
\begin{equation}\label{eq:2.22}
\left\|{\rm Op}_h(a_1){\rm Op}_h(a_2)-{\rm Op}_h(a_1a_2)\right\|_{L^2\to L^2}$$
 $$\lesssim \sum_{|\alpha_1|+|\beta_1|+|\beta_2|\le s'_d,\,|\alpha_1|+|\beta_1|\ge 1,\,1\le|\beta_2|\le \mu}h^{\frac{|\alpha_1|+|\beta_1|+|\beta_2|}{2}}C^{(1)}_{\alpha_1,\beta_1}C^{(2)}_{\beta_2}$$
 $$+h^{\mu/2}\sum_{|\alpha_1|+|\beta_1|\le s_d,\,|\beta_2|=\mu}
 h^{\frac{|\alpha_1|+|\beta_1|}{2}}C^{(1)}_{\alpha_1,\beta_1}C^{(2)}_{\beta_2}.
\end{equation}
\end{prop}

{\it Proof.} We will only prove (\ref{eq:2.19}) since (\ref{eq:2.20}) and (\ref{eq:2.22}) can be obtained in the same way. 
We approximate $a_1$ by the smooth function $a_{1,t}$ as in the proof of Proposition 2.4. Then the Taylor expansion allows us to write
the difference $a_1-a_{1,t}$ in the form
\begin{equation}\label{eq:2.23}
a_1-a_{1,t}=\sum_{s=1}^{\mu-1}t^sa_1^{(s)}+t^\mu a_{1,t}^{(\mu)},
\end{equation}
where
$$a_1^{(s)}=\sum_{|\beta|=s}c_\beta\partial_x^\beta a_1(x,\xi),\quad 1\le s\le\mu-1,$$
$$a_{1,t}^{(\mu)}=(\mu!)^{-1}\sum_{|\beta|=\mu}\int_{\mathbb{R}^{d-1}}\phi_0(w)w^\beta\partial_x^\beta a_1(x+t'w,\xi)dw.$$
Clearly, when $\mu=1$ the sum in the right-hand side of (\ref{eq:2.23}) is zero. We also have that $a_1^{(s)}$, $1\le s\le\mu-1$,
satisfy (\ref{eq:2.1}) with $k<0$ for all multi-indices
$\alpha$ and $\beta$ such that $|\beta|\le\mu-s$, while $a_{1,t}^{(\mu)}$ satisfies (\ref{eq:2.12}) with $k<0$ uniformly in $t$
for all multi-indices $\alpha$. In view of (\ref{eq:2.23}), we can write
\begin{equation}\label{eq:2.24}
{\rm Op}_h(a_1){\rm Op}_h(a_2)-{\rm Op}_h(a_1a_2)={\cal A}_0+\sum_{s=1}^{\mu-1}t^s{\cal A}_s+t^\mu{\cal A}_\mu,
\end{equation}
where 
$${\cal A}_0={\rm Op}_h(a_{1,t}){\rm Op}_h(a_2)-{\rm Op}_h(a_{1,t}a_2),$$
$${\cal A}_s={\rm Op}_h(a_1^{(s)}){\rm Op}_h(a_2)-{\rm Op}_h(a_1^{(s)}a_2),$$
$${\cal A}_\mu={\rm Op}_h(a_{1,t}^{(\mu)}){\rm Op}_h(a_2)-{\rm Op}_h(a_{1,t}^{(\mu)}a_2).$$
We take $t=h^{1/2}$ and use Proposition 2.2 to bound the norm of the operator ${\cal A}_0$. Taking into account the bounds
(\ref{eq:2.15}), we get
$$\left\|{\cal A}_0\right\|_{H_h^k\to L^2}$$
$$\lesssim \sum_{|\alpha_1|+|\beta_1|+|\alpha_2|+|\beta_2|\le s'_d,\,|\beta_1|\le\mu,\,|\alpha_1|+|\beta_1|\ge 1,\,|\alpha_2|+|\beta_2|\ge 1}h^{\frac{|\alpha_1|+|\beta_1|+|\alpha_2|+|\beta_2|}{2}}C^{(1)}_{\alpha_1,\beta_1}C^{(2)}_{\alpha_2,\beta_2}.$$
To bound the norm of the operator ${\cal A}_\mu$ we will use Propositions 2.1 and 2.3. We get
$$\left\|{\cal A}_\mu\right\|_{H_h^k\to L^2}\le \left\|{\rm Op}_h(a_{1,t}^{(\mu)})
\right\|_{L^2\to L^2}\left\|{\rm Op}_h(a_2)\right\|_{H_h^k\to L^2}
+\left\|{\rm Op}_h(a_{1,t}^{(\mu)}a_2)\right\|_{H_h^k\to L^2}$$
$$\lesssim\log(1+h^{-1})\sum_{|\alpha_1|\le d,\,|\beta_1|=\mu}C^{(1)}_{\alpha_1,\beta_1}
\sum_{|\alpha_2|+|\beta_2|\le s_d}h^{\frac{|\alpha_2|+|\beta_2|}{2}}C^{(2)}_{\alpha_2,\beta_2}$$ 
$$+\log(1+h^{-1})\sum_{|\alpha_1|+|\alpha_2|\le d,\,|\beta_1|=\mu,\,\beta_2=0}C^{(1)}_{\alpha_1,\beta_1}C^{(2)}_{\alpha_2,\beta_2}.$$
Now (\ref{eq:2.19}) follows from the above bounds and (\ref{eq:2.24}) by induction in $\mu$. 
\eproof

We will now use the above proposition to prove the following

\begin{prop} Let $a_1$ satisfy (\ref{eq:2.1}) with $k<0$ for all multi-indices
$\alpha$ and $\beta$ such that $|\beta|\le\mu_1$ with constants $C_{\alpha,\beta}^{(1)}$, and 
let $a_2$ satisfy (\ref{eq:2.1}) with $k<k'$ for all multi-indices
$\alpha$ and $\beta$ such that $|\beta|\le\mu_2$ with constants $C_{\alpha,\beta}^{(2)}$.  
 Then we have the bound
\begin{equation}\label{eq:2.25}
\left\|{\rm Op}_h(a_1){\rm Op}_h(a_2)-{\rm Op}_h(a_1a_2)\right\|_{H_h^{k'}\to L^2}$$
 $$\lesssim \sum_{|\alpha_1|+|\beta_1|+|\alpha_2|+|\beta_2|\le s'_d,\,|\beta_1|\le\mu_1,\,|\beta_2|\le\mu_2,\,|\alpha_1|+|\beta_1|\ge 1,\,|\alpha_2|+|\beta_2|\ge 1}h^{\frac{|\alpha_1|+|\beta_1|+|\alpha_2|+|\beta_2|}{2}}
 C^{(1)}_{\alpha_1,\beta_1}C^{(2)}_{\alpha_2,\beta_2}$$
 $$+h^{\mu_1/2}\log(1+h^{-1})\sum_{|\alpha_1|\le d,\,|\beta_1|=\mu_1,\,|\alpha_2|+|\beta_2|\le s_d,\,|\beta_2|\le\mu_2}h^{\frac{|\alpha_2|+|\beta_2|}{2}}C^{(1)}_{\alpha_1,\beta_1}C^{(2)}_{\alpha_2,\beta_2}$$ 
 $$+h^{\mu_1/2}\log(1+h^{-1})\sum_{|\alpha_1|+|\alpha_2|\le d,\,|\beta_1|=\mu_1,\,\beta_2=0}
 C^{(1)}_{\alpha_1,\beta_1}C^{(2)}_{\alpha_2,\beta_2}$$
 $$+h^{\mu_2/2}\log(1+h^{-1})\sum_{|\alpha_1|+|\beta_1|\le s_d,\,|\beta_1|\le\mu_1,\,|\alpha_2|\le d,\,|\beta_2|=\mu_2}h^{\frac{|\alpha_1|+|\beta_1|}{2}}C^{(1)}_{\alpha_1,\beta_1}C^{(2)}_{\alpha_2,\beta_2}$$
 $$+h^{\mu_2/2}\log(1+h^{-1})\sum_{|\alpha_1|+|\alpha_2|\le d,\,\beta_1=0,\,|\beta_2|=\mu_2}
 C^{(1)}_{\alpha_1,\beta_1}C^{(2)}_{\alpha_2,\beta_2}$$
 $$+h^{(\mu_1+\mu_2)/2}\log^2(1+h^{-1})\sum_{|\alpha_1|\le d,\,|\beta_1|=\mu_1,\,|\alpha_2|\le d,\,|\beta_2|=\mu_2}
 C^{(1)}_{\alpha_1,\beta_1}C^{(2)}_{\alpha_2,\beta_2}.
\end{equation}
Moreover, if $a_2$ is independent of the variable $\xi$ and satisfies (\ref{eq:2.21}) for all multi-indices $\beta$ such that $|\beta|\le\mu_2$, 
 then we have the bound
 \begin{equation}\label{eq:2.26}
\left\|{\rm Op}_h(a_1){\rm Op}_h(a_2)-{\rm Op}_h(a_1a_2)\right\|_{L^2\to L^2}$$
 $$\lesssim \sum_{|\alpha_1|+|\beta_1|+|\beta_2|\le s'_d,\,|\alpha_1|+|\beta_1|\ge 1,\,|\beta_1|\le\mu_1,\,1\le
 |\beta_2|\le\mu_2}h^{\frac{|\alpha_1|+|\beta_1|+|\beta_2|}{2}}
 C^{(1)}_{\alpha_1,\beta_1}C^{(2)}_{\beta_2}$$
 $$+h^{\mu_1/2}\log(1+h^{-1})\sum_{|\alpha_1|\le d,\,|\beta_1|=\mu_1}
 C^{(1)}_{\alpha_1,\beta_1}C^{(2)}_0$$
 $$+h^{\mu_2/2}\sum_{|\alpha_1|+|\beta_1|\le s_d,\,|\beta_1|\le\mu_1,\,|\beta_2|=\mu_2}h^{\frac{|\alpha_1|+|\beta_1|}{2}}C^{(1)}_{\alpha_1,\beta_1}C^{(2)}_{\beta_2}.
\end{equation}
\end{prop}

{\it Proof.} The proof is similar to the proof of (\ref{eq:2.19}) and we keep the same notations replacing $\mu$ by $\mu_1$. 
The only difference is that $a_2$ is no longer $C^\infty$ in $x$, so we have to bound the norms of the operators 
${\cal A}_0$ and ${\cal A}_{\mu_1}$ differently. To bound the norm of ${\cal A}_0$ we will use Proposition 2.5. By (\ref{eq:2.20}), we have
$$\left\|{\cal A}_0\right\|_{H_h^{k'}\to L^2}$$
$$\lesssim \sum_{|\alpha_1|+|\beta_1|+|\alpha_2|+|\beta_2|\le s'_d,\,|\beta_1|\le\mu_1\,|\beta_2|\le\mu_2,\,|\alpha_1|+|\beta_1|\ge 1,\,|\alpha_2|+|\beta_2|\ge 1}h^{\frac{|\alpha_1|+|\beta_1|+|\alpha_2|+|\beta_2|}{2}}
C^{(1)}_{\alpha_1,\beta_1}C^{(2)}_{\alpha_2,\beta_2}$$
 $$+h^{\mu_2/2}\log(1+h^{-1})\sum_{|\alpha_1|+|\beta_1|\le s_d,\,|\alpha_2|\le d,\,|\beta_1|\le\mu_1,\,|\beta_2|=\mu_2}h^{\frac{|\alpha_1|+|\beta_1|}{2}}C^{(1)}_{\alpha_1,\beta_1}C^{(2)}_{\alpha_2,\beta_2}$$
 $$+h^{\mu_2/2}\log(1+h^{-1})\sum_{|\alpha_1|+|\alpha_2|\le d,\,\beta_1=0,\,|\beta_2|=
 \mu_2}C^{(1)}_{\alpha_1,\beta_1}C^{(2)}_{\alpha_2,\beta_2}.$$
To bound the norm of the operator ${\cal A}_{\mu_1}$ we will use Propositions 2.3 and 2.4. We get
$$\left\|{\cal A}_{\mu_1}\right\|_{H_h^{k'}\to L^2}\le \left\|{\rm Op}_h(a_{1,t}^{(\mu_1)})
\right\|_{L^2\to L^2}\left\|{\rm Op}_h(a_2)\right\|_{H_h^{k'}\to L^2}
+\left\|{\rm Op}_h(a_{1,t}^{(\mu_1)}a_2)\right\|_{H_h^{k'}\to L^2}$$
$$\lesssim\log(1+h^{-1})\sum_{|\alpha_1|\le d,\,|\beta_1|=\mu_1}C^{(1)}_{\alpha_1,\beta_1}
\sum_{|\alpha_2|+|\beta_2|\le s_d,\,|\beta_2|\le\mu_2}h^{\frac{|\alpha_2|+|\beta_2|}{2}}C^{(2)}_{\alpha_2,\beta_2}$$ 
$$+h^{\mu_2/2}\log^2(1+h^{-1})\sum_{|\alpha_1|\le d,\,|\beta_1|=\mu_1}C^{(1)}_{\alpha_1,\beta_1}
\sum_{|\alpha_2|\le d,\,|\beta_2|=\mu_2}C^{(2)}_{\alpha_2,\beta_2}$$ 
$$+\log(1+h^{-1})\sum_{|\alpha_1|+|\alpha_2|\le d,\,|\beta_1|=\mu_1,\,\beta_2=0}C^{(1)}_{\alpha_1,\beta_1}C^{(2)}_{\alpha_2,\beta_2}.$$
Now (\ref{eq:2.25}) follows from the above bounds and (\ref{eq:2.24}) by induction in $\mu_1$. The bound (\ref{eq:2.26})
can be derived from (\ref{eq:2.22}) in the same way.
\eproof

It is well-known that one can define $h-\Psi$DOs with $C^\infty$ smooth symbols on an arbitrary compact manifold without boundary
using the definition on the Euclidean space (e.g. see Chapter 18 of \cite{kn:H}). 
Roughly speaking, an $h-\Psi$DO on a manifold of dimension $d-1$, say $\Gamma$, is a finite sum of 
$h-\Psi$DOs on $\mathbb{R}^{d-1}$ with symbols compactly supported in $x$. In other words, studying $h-\Psi$DOs on compact manifolds is
reduced to studying $h-\Psi$DOs on the Euclidean space with symbols compactly supported in $x$. In the same way, we can define 
$h-\Psi$DOs on $\Gamma$ with symbols of low regularity in $x$. To be more precise, 
we cover $\Gamma$ by a finite number of open sets $U_j$, $j=1,...,J$, such that
$$\kappa_j:U_j\to \widetilde U_j\subset \mathbb{R}^{d-1}$$
are diffeomorphisms. Then we can associate to $\kappa_j$ a diffeomorphism 
$${\cal K}_j:T^*U_j\to T^*\widetilde U_j\subset T^*\mathbb{R}^{d-1}$$
such that
$${\cal K}_j(x,\xi)=(\kappa_j(x),\,^t\kappa'_j(x)\xi).$$
Let $\psi_j\in C^\infty(U_j)$ be such that $\sum_{j=1}^J\psi_j=1$. If $a$ is a function on $T^*\Gamma$, 
we can write it in the form $a=\sum_{j=1}^Ja_j$, where $a_j=a\psi_j$. 
Define the function $a_j\circ{\cal K}_j^{-1}$ on $T^*\widetilde U_j$
 by 
 $$(a_j\circ {\cal K}_j^{-1})(y,\eta)=a_j({\cal K}_j^{-1}(y,\eta)).$$ 
 Then we define the operators ${\rm Op}_h(a_j\circ{\cal K}_j^{-1})$
 as in the begining of the section. We now define ${\rm Op}_h(a)$ in terms of these operators as follows
 $${\rm Op}_h(a)f=\sum_{j=1}^J\sum_{\ell=1}^J{\rm Op}_h(\psi_ja)\psi_\ell f
 =\sum_{j=1}^J\sum_{\ell=1}^J\left({\rm Op}_h(a_j\circ{\cal K}_j^{-1})((\psi_\ell f)\circ\kappa_\ell^{-1})\right)\circ\kappa_j.$$
 We will say that $a\in S^k(\Gamma)$ if $a_j\circ{\cal K}_j^{-1}\in S^k$ for all $j$. Similarly, we can extend the other classes 
 of symbols on $T^*\mathbb{R}^{d-1}$ above 
 to symbols on $T^*\Gamma$. It is then clear that the above propositions extend to $h-\Psi$DOs on $\Gamma$ with the spaces
 $L^2$ and $H_h^k$ replaced by $L^2(\Gamma)$ and $H_h^k(\Gamma)$, respectively. To simplify the notations in
 the formulas that follow,
 we will omit the diffeomorphisms and will identify $T^*U_j$ with $T^*\widetilde U_j$.

\section{A priori estimates}

In this section we will prove a priori estimates for the solution to the equation
\begin{equation}\label{eq:3.1}
\left\{
\begin{array}{l}
(h^2\nabla c(x)\nabla+zn(x))u=hv\quad \mbox{in}\quad\Omega,\\
u=0\quad\mbox{on}\quad\Gamma,\\
\end{array}
\right.
\end{equation}
where $c,n\in L^\infty(\Omega)$ are real-valued functions satisfying $c(x)\ge b_0$, $n(x)\ge b_0$ for some constant $b_0>0$,
$0<h\ll 1$ is a semiclassical parameter and $z\in Z^+\cup Z^-$, where 
$$Z^\pm=\{z\in \mathbb{C}:|z|=1,\,\pm{\rm Re}\,z\ge 0\}.$$
Set $\theta=|{\rm Im}\,z|$ if $z\in Z^+$ and $\theta=1$ if $z\in Z^-$. 
We will prove the following

\begin{Theorem} \label{3.1} Suppose that $c, n\in C^1(\overline\Omega_\delta)$ for some $0<\delta\ll 1$. Let $\theta\ge h$ and let $u\in H^2(\Omega)$ satisfy equation (\ref{eq:3.1}). 
Then the function $g=h\partial_\nu u|_\Gamma$ satisfies the estimate
\begin{equation}\label{eq:3.2}
\|g\|_{L^2(\Gamma)}\lesssim h^{1/2}\theta^{-1/2}\|v\|_{L^2(\Omega)}.
\end{equation}
\end{Theorem}

{\it Proof.} We will first prove the following

\begin{lemma} \label{3.2} We have the estimate
\begin{equation}\label{eq:3.3}
\|u\|_{H_h^1(\Omega)}\lesssim h\theta^{-1}\|v\|_{L^2(\Omega)}.
\end{equation}
\end{lemma}

{\it Proof.} By the Green formula we have
\begin{equation}\label{eq:3.4}
\left\langle znu-hv,u\right\rangle_{L^2(\Omega)}=
\left\langle -h^2\nabla c\nabla u,u\right\rangle_{L^2(\Omega)}=\int_\Omega c|h\nabla u|^2.
\end{equation}
 Taking the imaginary part we get the identity
$${\rm Im}\,z\left\|n^{1/2}u\right\|^2_{L^2(\Omega)}={\rm Im}\,\left\langle hv,u\right\rangle_{L^2(\Omega)},$$
which implies
\begin{equation}\label{eq:3.5}
\|u\|_{L^2(\Omega)}\lesssim h|{\rm Im}\,z|^{-1}\|v\|_{L^2(\Omega)}.
\end{equation}
Taking the real part of (\ref{eq:3.4}) we get
\begin{equation}\label{eq:3.6}
\int_\Omega c|h\nabla u|^2
={\rm Re}\,z\left\langle nu,u\right\rangle_{L^2(\Omega)}-{\rm Re}\,\left\langle hv,u\right\rangle_{L^2(\Omega)}$$
$$\le ({\rm Re}\,z+\varepsilon)\left\langle nu,u\right\rangle_{L^2(\Omega)}+{\cal O}_\varepsilon(h^2)\|v\|^2_{L^2(\Omega)}
\end{equation}
for every $0<\varepsilon \le 1$. In particular, (\ref{eq:3.6}) implies
\begin{equation}\label{eq:3.7}
\int_\Omega |h\nabla u|^2\lesssim\|u\|^2_{L^2(\Omega)}+ h^2\|v\|^2_{L^2(\Omega)}.
\end{equation}
By (\ref{eq:3.5}) and (\ref{eq:3.7}) we obtain 
\begin{equation}\label{eq:3.8}
\|u\|_{H_h^1(\Omega)}\lesssim h|{\rm Im}\,z|^{-1}\|v\|_{L^2(\Omega)}.
\end{equation}
When $|{\rm Re}\,z|\le 1/2$ we have $|{\rm Im}\,z|\ge 1/2$, so in this case 
 (\ref{eq:3.3}) follows from (\ref{eq:3.8}). When ${\rm Re}\,z\ge 1/2$ we have $|{\rm Im}\,z|=\theta$, so in this case 
 (\ref{eq:3.3}) again follows from (\ref{eq:3.8}). When ${\rm Re}\,z\le -1/2$ we have $\theta=1$, so
 in this case (\ref{eq:3.3}) follows from (\ref{eq:3.6}) 
 with $\varepsilon=1/4$.
\eproof

Let ${\cal V}\subset\mathbb{R}^d$ be a small open domain such that ${\cal V}^0:={\cal V}\cap\Gamma\neq\emptyset$. 
Let $(x_1,x')\in {\cal V}^+:={\cal V}\cap\Omega$, $0<x_1\ll 1$, $x'=(x_2,...,x_d)\in{\cal V}^0$, be the local normal geodesic coordinates near the boundary. In these coordinates the principal symbol of the operator $-\Delta$ is equal to 
$\xi_1^2+r(x,\xi')$, where $(\xi_1,\xi')$ are the dual variables to $(x_1,x')$, and $r$ is a homogeneous polynomial of order two and satisfies $C_1|\xi'|^2\le r\le C_2|\xi'|^2$
with some constants $C_1,C_2>0$. 
Therefore, the principal symbol of the positive Laplace-Beltrami operator on $\Gamma$ is equal to
$r_0(x',\xi')=r(0,x',\xi')$. 

Let ${\cal V}_1\subset{\cal V}$ be a small open domain such that ${\cal V}_1^0:={\cal V_1}\cap\Gamma\neq\emptyset$. 
Choose a function $\psi\in C_0^\infty({\cal V})$, $0\le\psi\le 1$, such that $\psi=1$ on ${\cal V}_1$. 
 Then the function $u^\flat:=\psi u$ satisfies the equation 
\begin{equation}\label{eq:3.9}
\left\{
\begin{array}{l}
(h^2\nabla c(x)\nabla+zn)u^\flat=hv^\flat\quad \mbox{in}\quad\Omega,\\
u^\flat=0\quad\mbox{on}\quad\Gamma,\\
\end{array}
\right.
\end{equation}
where $v^\flat=\psi v+h[\nabla c\nabla,\psi]u$ satisfies
\begin{equation}\label{eq:3.10}
\|v^\flat\|_{L^2(\Omega)}\lesssim\|v\|_{L^2(\Omega)}+\|u\|_{H^1_h(\Omega)}.
\end{equation}
 We will now write the operator $-h^2\nabla c(x)\nabla$ in the coordinates $x=(x_1,x')$.  
 Denote ${\cal D}_{x_j}=-ih\partial_{x_j}$. We can write
\begin{equation}\label{eq:3.11}
-h^2\nabla c(x)\nabla=c(x){\cal D}_{x_1}^2+c(x)r(x,{\cal D}_{x'})+h{\cal R}(x,{\cal D}_{x}),
\end{equation}
where ${\cal R}$ is a first-order differential operator with $L^\infty$ coefficients. Denote
$$\langle f,g\rangle_0=\int f\overline g dx',\quad \|f\|_0^2=\int|f|^2dx',$$
and introduce the function
$$F(x_1)=\left\|{\cal D}_{x_1}u^\flat\right\|_0^2
-\left\langle r(x_1,\cdot,{\cal D}_{x'})u^\flat,u^\flat\right\rangle_0+{\rm Re}\,z\left\langle \widetilde n(x_1,\cdot)u^\flat,u^\flat\right\rangle_0,$$
where $\widetilde n=c^{-1}n\in C^1(\overline\Omega_\delta)$. 
Since $u^\flat|_{x_1=0}=0$, we have 
\begin{equation}\label{eq:3.12}
F(0)=\left\|{\cal D}_{x_1}u^\flat|_{x_1=0}\right\|_0^2.
\end{equation}
 On the other hand,
\begin{equation}\label{eq:3.13}
F(0)=-\int_0^{\delta}F'(x_1)dx_1
\end{equation}
for some constant $\delta>0$, 
where $F'$ denotes the first derivative with respect to $x_1$. We will now bound $F(0)$ from above. 
To this end we will compute $F'(x_1)$ using that $u^\flat$ satisfies (\ref{eq:3.9}) together with (\ref{eq:3.11}).
We have
$$F'(x_1)=-2{\rm Re}\,\left\langle ({\cal D}_{x_1}^2+r-{\rm Re}\,z\widetilde n)u^\flat,
\partial_{x_1}u^\flat\right\rangle_0-\left\langle (r'-{\rm Re}\,z\widetilde n')u^\flat,u^\flat\right\rangle_0$$
$$=2h^{-1}{\rm Im}\,\left\langle (-h^2\nabla c(x)\nabla+{\rm Re}\,zn-h{\cal R})u^\flat,c^{-1}{\cal D}_{x_1}u^\flat\right\rangle_0$$
$$-\left\langle (r'-{\rm Re}\,z\widetilde n')u^\flat,u^\flat\right\rangle_0$$
$$=2{\rm Im}\,\left\langle (v^\flat-ih^{-1}{\rm Im}\,znu^\flat-{\cal R}u^\flat),c^{-1}{\cal D}_{x_1}u^\flat\right\rangle_0$$
$$-\left\langle (r'-{\rm Re}\,z\widetilde n')u^\flat,u^\flat\right\rangle_0.$$
Hence
$$|F'(x_1)|\lesssim h\theta^{-1}\|v^\flat\|_0^2+\theta h^{-1}\sum_{\ell=0}^1\|{\cal D}_{x_1}^\ell u^\flat\|_0^2+\sum_{|\alpha|\le 1}\|{\cal D}_{x}^\alpha  u^\flat\|_0^2.$$
Using this estimate together with (\ref{eq:3.10}), (\ref{eq:3.13}) and Lemma 3.2 we obtain
\begin{equation}\label{eq:3.14}
F(0)\le \int_0^{\delta}|F'(x_1)|dx_1\lesssim h\theta^{-1}\|v\|_{L^2(\Omega)}^2
+(1+\theta h^{-1})\|u\|_{H_h^1(\Omega)}^2$$ 
$$\lesssim (h\theta^{-1}+h^2\theta^{-2})\|v\|_{L^2(\Omega)}^2\lesssim h\theta^{-1}\|v\|_{L^2(\Omega)}^2.
\end{equation}
Observe now that
 $${\cal D}_{x_1}u^\flat|_{x_1=0}=\psi_0{\cal D}_{x_1}u|_{x_1=0},\quad {\cal D}_{x'}u|_{x_1=0}=0,$$
 where $\psi_0=\psi|_{x_1=0}$ is supported in ${\cal V}^0$ and such that $\psi_0=1$ on ${\cal V}_1^0$. 
 Therefore, by (\ref{eq:3.12}) and (\ref{eq:3.14}),
$$
\left\|\psi_0{\cal D}_{x_1}u|_{x_1=0}\right\|_0\lesssim h^{1/2}\theta^{-1/2}\|v\|_{L^2(\Omega)},
$$
which clearly implies 
\begin{equation}\label{eq:3.15}
\left\|\psi_0g\right\|_0\lesssim h^{1/2}\theta^{-1/2}\|v\|_{L^2(\Omega)}.
\end{equation}
Since $\Gamma$ is compact, there exist a finite number of smooth functions $\psi_i$, $0\le\psi_i\le 1$, $i=1,...,I,$ such that 
$1=\sum_{i=1}^I\psi_i$ and (\ref{eq:3.15}) holds with $\psi_0$ replaced by each $\psi_i$. Therefore, 
 the estimate (\ref{eq:3.2}) is obtained by summing up all such estimates (\ref{eq:3.15}).
\eproof

\section{Approximation of the Dirichlet-to-Neumann map}

Given $f\in H^1(\Gamma)$ let $u$ solve the equation 
\begin{equation}\label{eq:4.1}
\left\{
\begin{array}{l}
(h^2\nabla c(x)\nabla+zn(x))u=0\quad \mbox{in}\quad\Omega,\\
u=f\quad\mbox{on}\quad\Gamma,\\
\end{array}
\right.
\end{equation}
where $c,n,h,z$ are as in previous section. We define the semiclassical Dirichlet-to-Neumann map 
$${\cal N}(h,z):H^1(\Gamma)\to L^2(\Gamma)$$
by
$${\cal N}(h,z)f:=-ih\partial_\nu u|_\Gamma.$$
We would like to approximate ${\cal N}(h,z)$ by an $h-\Psi$DO on $\Gamma$ similarly to the $C^\infty$ case. To this end,
introduce the function
$$\rho(x',\xi',z)=\sqrt{-r_0(x',\xi')+z\widetilde n_0(x')},\quad {\rm Im}\,\rho>0,$$
where $\widetilde n_0=\widetilde n|_\Gamma$, $\widetilde n:=c^{-1}n$, and $r_0$ is the principal symbol of the Laplace-Beltrami operator on $\Gamma$ written
in the coordinates $(x',\xi')\in T^*\Gamma$. 
Let $\eta\in C^\infty(T^*\Gamma)$ be such that $\eta=1$ for $r_0\le C_0$, $\eta=0$ for $r_0\ge 2C_0$, where $C_0>0$ does not depend on $h$. 
It is easy to see (e.g. see Lemma 3.1 of \cite{kn:V1}) that taking $C_0$ big enough we can arrange 
\begin{equation}\label{eq:4.2}
C_1\theta^{1/2}\le |\rho|\le C_2,\quad {\rm Im}\,\rho\ge C_3|\theta||\rho|^{-1}\ge C_4|\theta|
\end{equation}
 for $(x',\xi')\in{\rm supp}\,\eta$, and 
 \begin{equation}\label{eq:4.3}
 |\rho|\ge {\rm Im}\,\rho\ge C_5|\xi'|
 \end{equation}
  for $(x',\xi')\in{\rm supp}\,(1-\eta)$ with some constants $C_j>0$. In other words, ${\rm supp}\,(1-\eta)$ is contained in the
  elliptic region of the boundary value problem (\ref{eq:4.1}). 
  Our goal is to prove the following

\begin{Theorem} Let $\theta\ge h^{2/5}$. Suppose that $c,n\in C^1(\Omega_\delta)$ and 
$\widetilde n_0\in C^\mu(\Gamma)$ with an integer $\mu\ge 2$.
 Then for every $f\in H^1(\Gamma)$ 
we have the estimate
\begin{equation}\label{eq:4.4}
\left\|{\cal N}(z,h)f-{\rm Op}_h(\rho)f\right\|_{L^2(\Gamma)}\lesssim h^{3/5}\|f\|_{H_h^1(\Gamma)}$$  $$+h\theta^{-5/2}\left(1+h^{\mu/2}\log(h^{-1})\theta^{-d-\mu}\right)\|f\|_{H_h^1(\Gamma)}. 
\end{equation}
Suppose in addition that $c\equiv 1$. Then we have the better estimate 
\begin{equation}\label{eq:4.5}
\left\|{\cal N}(z,h)f-{\rm Op}_h(\rho+hq)f\right\|_{L^2(\Gamma)}\lesssim h\theta^{-5/2}\left(1+h^{\mu/2}\log(h^{-1})\theta^{-d-\mu}\right)\|f\|_{H_h^{-1}(\Gamma)},  
\end{equation}
where $q\in C^\infty(T^*\Gamma)$ is independent of $z$, $h$ and the function $n$.
 \end{Theorem}

{\it Proof.} We will build a parametrix for the solutions of the equation (\ref{eq:4.1}). This parametrix will not be as good as
that one built in the $C^\infty$ case, but will suffice for the proof of the above estimates. 
  Let $(x_1,x')\in {\cal V}^+$ be the local normal geodesic coordinates near the boundary. 
Take a function $\chi\in C^\infty(T^*\Gamma)$, $0\le\chi\le 1$, such that $\pi_{x'}({\rm supp}\,\chi)\subset {\cal V}^0$, where 
$\pi_{x'}:T^*\Gamma\to\Gamma$ denotes the projection $(x',\xi')\to x'$. Moreover, we require that either $\chi$ is of compact support or
 $\chi\in S^0(\Gamma)$ with ${\rm supp}\,\chi\subset{\rm supp}(1-\eta)$. 
We will be looking for a parametrix of the solution to 
equation (\ref{eq:4.1}) in the form
$$\widetilde u=\phi(x_1/\delta)(2\pi h)^{-d+1}\int\int e^{\frac{i}{h}(\langle y',\xi'\rangle+\varphi(x,\xi',z))}a(x,\xi',z)f(y')d\xi'dy'$$
where $\phi\in C_0^\infty(\mathbb{R})$,
$\phi(t)=1$ for $|t|\le 1/2$, $\phi(t)=0$ for $|t|\ge 1$, $0<\delta\ll 1$ being a parameter independent of $h$ and $z$.  
 We require that $\widetilde u$ satisfies the boundary condition
$\widetilde u={\rm Op}_h(\chi)f$ on $x_1=0$. We take $a=\chi(x',\xi')$. Furhtermore, we choose the phase function in the form
$$\varphi=-\langle x',\xi'\rangle+x_1\rho.$$
Then it is easy to see that $\varphi$ satisfies the following eikonal equation
\begin{equation}\label{eq:4.6}
(\partial_{x_1}\varphi)^2+r(x,\nabla_{x'}\varphi)-z\widetilde n(x)=x_1\Phi,
\end{equation}
where $|\Phi|$ is bounded as $x_1\to 0$. By assumption, we have 
$$\widetilde n(x)=\widetilde n_0(x')+x_1\widetilde n^\sharp(x),\quad \widetilde n_0\in C^\mu,\,\widetilde n^\sharp\in L^\infty.$$
Furthermore, it is clear that the function $r$ can be written in the form
$$r(x,\xi')=\langle R(x)\xi',\xi'\rangle$$
where $R$ is a $C^\infty$ smooth $(d-1)\times(d-1)$ matrix-valued function. Hence we can write
$$R(x)=R_0(x')+x_1R^\sharp(x)$$
where $R_0=R|_{x_1=0}$ and $R^\sharp$ are $C^\infty$ smooth functions. Thus we obtain
$$\Phi=-z\widetilde n^\sharp+\langle R^\sharp(x)\xi',\xi'\rangle-2\langle(R_0(x')+x_1R^\sharp(x))\xi',\nabla_{x'}\rho\rangle 
+x_1\langle R(x)\nabla_{x'}\rho,\nabla_{x'}\rho\rangle.$$
Observe now that 
\begin{equation}\label{eq:4.7} 
\nabla_{x'}\rho=(2\rho)^{-1}\left(-\nabla_{x'}r_0+z\nabla_{x'}\widetilde n_0\right).
\end{equation}
Hence we can write the function $\Phi$ in the form
$$\Phi=-z\widetilde n^\sharp+\langle R^\sharp(x)\xi',\xi'\rangle
+\rho^{-1}\langle (R_0(x')+x_1R^\sharp(x))\xi',\nabla_{x'}r_0-z\nabla_{x'}\widetilde n_0\rangle$$
$$+x_1(2\rho)^{-2}\langle R(x)(\nabla_{x'}r_0-z\nabla_{x'}\widetilde n_0),\nabla_{x'}r_0-z\nabla_{x'}\widetilde n_0\rangle.$$
Observe now that $\widetilde u$ satisfies the equation 
$$(h^2\nabla c(x)\nabla+zn(x))\widetilde u=h\widetilde v,$$
where the function $\widetilde v$ is of the form
$$\widetilde v=(2\pi h)^{-d+1}\int\int e^{\frac{i}{h}\langle y'-x',\xi'\rangle}
e^{ix_1\rho/h}A(x,\xi',h,z)f(y')d\xi'dy'$$ 
$$={\rm Op}_h\left(e^{ix_1\rho/h}A\right)f$$
where
$$A=h^{-1}e^{-i\varphi/h}(h^2\nabla c(x)\nabla+zn(x))\left(\phi(x_1/\delta)e^{i\varphi/h}a\right).$$ 
To compute the function $A$ observe that the operator $-\nabla c\nabla$ can be written in the coordinates $(x_1,x')$ in the form
$$-\nabla c\nabla=cD_{x_1}^2+c\left\langle R(x)D_{x'},D_{x'}\right\rangle+
\left\langle Q(x),D_{x}\right\rangle,$$
where $D_{x_1}=-i\partial_{x_1}$, $D_{x'}=-i\nabla_{x'}$, $D_{x}=-i\nabla_{x}=(D_{x_1},D_{x'})$, 
$R$ is the smooth matrix-valued function as above, and $Q=(Q_1,...,Q_d)$ with scalar-valued functions
$Q_j\in L^\infty$. Denote $\widetilde Q=(Q_2,...,Q_d)$. Hence we can write
$$-h^2\nabla c\nabla=c{\cal D}_{x_1}^2+hQ_1(x){\cal D}_{x_1}+c\left\langle R(x){\cal D}_{x'},{\cal D}_{x'}\right\rangle+
h\left\langle \widetilde Q(x),{\cal D}_{x'}\right\rangle,$$
where ${\cal D}_{x_1}=hD_{x_1}$, ${\cal D}_{x'}=hD_{x'}$. Therefore, the function $A$ can be decomposed as $A_1+\phi(x_1/\delta)A_2$,
where
$$A_1=e^{-i\varphi/h}\left[ch^{-1}{\cal D}_{x_1}^2+Q_1(x){\cal D}_{x_1},\phi(x_1/\delta)\right]\left(e^{i\varphi/h}a\right)$$
$$=-2ic\delta^{-1}\phi'(x_1/\delta)e^{-i\varphi/h}{\cal D}_{x_1}\left(e^{i\varphi/h}a\right)$$
$$-hc\delta^{-2}\phi''(x_1/\delta)a-ihQ_1\delta^{-1}\phi'(x_1/\delta)a$$
$$=-2ic\delta^{-1}\phi'(x_1/\delta)\left(a\partial_{x_1}\varphi+{\cal D}_{x_1}a\right)$$
$$-hc\delta^{-2}\phi''(x_1/\delta)a-ihQ_1\delta^{-1}\phi'(x_1/\delta)a$$
$$=-2ic\delta^{-1}\phi'(x_1/\delta)\chi\rho-hc\delta^{-2}\phi''(x_1/\delta)\chi-ihQ_1\delta^{-1}\phi'(x_1/\delta)\chi,$$
and 
$$A_2=-h^{-1}cx_1\Phi-2ic\left(\partial_{x_1}\varphi\partial_{x_1}a+\langle\nabla_{x'}\varphi,\nabla_{x'}a\rangle\right)$$
$$+ica\left(D_{x_1}^2+\left\langle R(x)D_{x'},D_{x'}\right\rangle\right)\varphi
+hc\left(D_{x_1}^2+\left\langle R(x)D_{x'},D_{x'}\right\rangle\right)a$$
$$+Q_1\left(a\partial_{x_1}\varphi+{\cal D}_{x_1}a\right)+\left\langle\widetilde Q,\nabla_{x'}\varphi\right\rangle a
+\left\langle\widetilde Q,{\cal D}_{x'}a\right\rangle$$
$$=-h^{-1}cx_1\Phi-2ic\langle-\xi'+x_1\nabla_{x'}\rho,\nabla_{x'}\chi\rangle-icx_1\chi\left\langle R(x)\nabla_{x'},\nabla_{x'}\right\rangle\rho
$$
$$-hc\left\langle R(x)\nabla_{x'},\nabla_{x'}\right\rangle\chi+Q_1\chi\rho+\left\langle\widetilde Q,-\xi'+x_1\nabla_{x'}\rho\right\rangle\chi
-ih\left\langle\widetilde Q,\nabla_{x'}\chi\right\rangle$$
$$=-h^{-1}cx_1\Phi+A_3.$$
In view of (\ref{eq:4.7}), the function $A_3$ can be written in the form
$$A_3=2ic\langle\xi',\nabla_{x'}\chi\rangle+ix_1\rho^{-1}c\langle\nabla_{x'}r_0-z\nabla_{x'}\widetilde n_0,\nabla_{x'}\chi\rangle$$
$$-hc\left\langle R(x)\nabla_{x'},\nabla_{x'}\right\rangle\chi+Q_1\chi\rho-\left\langle\widetilde Q,\xi'\right\rangle\chi
-ih\left\langle\widetilde Q,\nabla_{x'}\chi\right\rangle$$
$$-x_1(2\rho)^{-1}\left\langle\widetilde Q,\nabla_{x'}r_0-z\nabla_{x'}\widetilde n_0\right\rangle\chi-icx_1\chi\left\langle R(x)\nabla_{x'},\nabla_{x'}\right\rangle\rho.$$
The last term can also be expressed in terms of negative powers of $\rho$. Indeed, if $2\le\ell,j\le d$, we have 
$$\partial_{x_\ell}\partial_{x_j}\rho=\partial_{x_\ell}\left((2\rho)^{-1}\left(-\partial_{x_j}r_0+z\partial_{x_j}\widetilde n_0\right)\right)$$
$$=(2\rho)^{-1}\left(-\partial_{x_\ell}\partial_{x_j}r_0+z\partial_{x_\ell}\partial_{x_j}\widetilde n_0\right)+
2^{-1}\rho^{-2}\left(\partial_{x_j}r_0-z\partial_{x_j}\widetilde n_0\right)\partial_{x_\ell}\rho$$
$$=(2\rho)^{-1}\left(-\partial_{x_\ell}\partial_{x_j}r_0+z\partial_{x_\ell}\partial_{x_j}\widetilde n_0\right)-
2^{-2}\rho^{-3}\left(\partial_{x_j}r_0-z\partial_{x_j}\widetilde n_0\right)\left(\partial_{x_\ell}r_0-z\partial_{x_\ell}\widetilde n_0\right).$$
It is easy to see from the above expressions that we have the following 

\begin{lemma} The functions $\Phi$ and $A_3$ are of the form
\begin{equation}\label{eq:4.8} 
\Phi=b_1+b_2\rho^{-1}+x_1\rho^{-1}\left(b_3+b_4\rho^{-1}\right),
\end{equation}
\begin{equation}\label{eq:4.9} 
A_3=b_5+hb_6+b_7\rho+x_1\rho^{-1}\left(b_8+b_9\rho^{-1}+b_{10}\rho^{-2}\right),
\end{equation}
where $b_j$ are functions on $T^*\Gamma$ independent of $h$, $\rho$ and depending on $x_1$ and $z$. Moreover,
 each function $b_j$ is of the form 
\begin{equation}\label{eq:4.10} 
b_j=\sum_{|\alpha|\le\ell_j}\omega_{j,\alpha}(x,z)\chi_{j,\alpha}(x',\xi')\xi'^\alpha,
\end{equation}
where $\chi_{j,\alpha}$ are either $\chi$ or derivatives of $\chi$, 
$\omega_{j,\alpha}\in L^\infty$ uniformly in $z$, and $0\le \ell_j\le 4$. More precisely, we have
$\ell_1=2$, $\ell_2=\ell_3=3$, $\ell_4=4$, $\ell_5=1$, $\ell_6=\ell_7=0$, $\ell_8=2$, $\ell_9=3$, $\ell_{10}=4$. 
\end{lemma}

Using this lemma we will prove the following

\begin{prop} Let $0<x_1\le \delta$. If $\chi$ is of compact support, there exists a constant $C>0$ such that, 
given any $N\ge 0$ independent of $h$, we have the estimate
\begin{equation}\label{eq:4.11}
\left\|\left({\rm Op}_h\left(e^{ix_1\rho/h}A_2\right)f\right)(x_1,\cdot)\right\|_{L^2(\Gamma)}$$
$$\lesssim \left(1+h^{\mu/2}\log(h^{-1})\theta^{-d-\mu}\right)\theta^{-3/2}e^{-Cx_1\theta/2h}\|f\|_{H_h^{-N}(\Gamma)}.
\end{equation}
If ${\rm supp}\,\chi\subset{\rm supp}(1-\eta)$, given any $0<\epsilon<1/2$ independent of $h$, we have the estimate
\begin{equation}\label{eq:4.12}
\left\|\left({\rm Op}_h\left(e^{ix_1\rho/h}A_2\right)f\right)(x_1,\cdot)\right\|_{L^2(\Gamma)}\lesssim h^{(1-3\epsilon)/2}x_1^{-1/2+\epsilon}\|f\|_{H_h^1(\Gamma)}.
\end{equation}
\end{prop}

{\it Proof.} Clearly, $\widetilde\chi_{j,\alpha}:=\chi_{j,\alpha}\xi'^\alpha\in C^\infty(T^*\Gamma)$ is compactly supported
if so is $\chi$, and $\widetilde\chi_{j,\alpha}\in S^{|\alpha|}(\Gamma)$, ${\rm supp}\,\widetilde\chi_{j,\alpha}\subset{\rm supp}(1-\eta)$
if $\chi\in S^0(\Gamma)$, ${\rm supp}\,\chi\subset{\rm supp}(1-\eta)$. 
In view of Lemma 4.2 we can write
$${\rm Op}_h\left(e^{ix_1\rho/h}A_2\right)=
-h^{-1}x_1\sum_{|\alpha|\le 2}c\omega_{1,\alpha}{\rm Op}_h\left(e^{ix_1\rho/h}\widetilde\chi_{1,\alpha}\right)$$
$$-h^{-1}x_1\sum_{|\alpha|\le 3}c\omega_{2,\alpha}{\rm Op}_h\left(e^{ix_1\rho/h}\rho^{-1}\widetilde\chi_{2,\alpha}\right)$$
$$-h^{-1}x_1^2\sum_{|\alpha|\le 3}c\omega_{3,\alpha}{\rm Op}_h\left(e^{ix_1\rho/h}\rho^{-1}\widetilde\chi_{3,\alpha}\right)$$ 
$$-h^{-1}x_1^2\sum_{|\alpha|\le 4}c\omega_{4,\alpha}{\rm Op}_h\left(e^{ix_1\rho/h}\rho^{-2}\widetilde\chi_{4,\alpha}\right)$$
$$+\sum_{|\alpha|\le 1}\omega_{5,\alpha}{\rm Op}_h\left(e^{ix_1\rho/h}\widetilde\chi_{5,\alpha}\right)
+h\omega_{6,0}{\rm Op}_h\left(e^{ix_1\rho/h}\widetilde\chi_{6,0}\right)
+\omega_{7,0}{\rm Op}_h\left(e^{ix_1\rho/h}\rho\widetilde\chi_{7,0}\right)$$
$$+x_1\sum_{|\alpha|\le 2}\omega_{8,\alpha}{\rm Op}_h\left(e^{ix_1\rho/h}\rho^{-1}\widetilde\chi_{8,\alpha}\right)$$ 
$$+x_1\sum_{|\alpha|\le 3}\omega_{9,\alpha}{\rm Op}_h\left(e^{ix_1\rho/h}\rho^{-2}\widetilde\chi_{9,\alpha}\right)$$ 
$$+x_1\sum_{|\alpha|\le 4}\omega_{10,\alpha}{\rm Op}_h\left(e^{ix_1\rho/h}\rho^{-3}\widetilde\chi_{10,\alpha}\right).$$
Hence
$$\left\|{\rm Op}_h\left(e^{ix_1\rho/h}A_2\right)f\right\|\lesssim 
\sum_{\ell=1}^2\sum_{k=\ell-1}^\ell h^{-1}x_1^\ell\sum_{|\alpha|\le k+2}
\left\|{\rm Op}_h\left(e^{ix_1\rho/h}\rho^{-k}\widetilde\chi_{j(k,\ell),\alpha}\right)f\right\|$$
$$+\sum_{|\alpha|\le 1}\left\|{\rm Op}_h\left(e^{ix_1\rho/h}\widetilde\chi_{5,\alpha}\right)f\right\|
+h\left\|{\rm Op}_h\left(e^{ix_1\rho/h}\widetilde\chi_{6,0}\right)f\right\|
+\left\|{\rm Op}_h\left(e^{ix_1\rho/h}\rho\widetilde\chi_{7,0}\right)f\right\|$$
$$+x_1\sum_{k=1}^3\sum_{|\alpha|\le k+1}
\left\|{\rm Op}_h\left(e^{ix_1\rho/h}\rho^{-k}\widetilde\chi_{7+k,\alpha}\right)f\right\|.$$
It is easy to see now that the proposition follows from

\begin{prop} Let $k\ge 0$ be an integer or $k=-1$ and let $\ell\ge 0$. If $\widetilde\chi\in C^\infty(T^*\Gamma)$ is of compact support, there exists a constant $C>0$ such that, 
given any $N\ge 0$ independent of $h$, we have the bound
\begin{equation}\label{eq:4.13}
\left\|{\rm Op}_h\left(e^{ix_1\rho/h}\rho^{-k}\widetilde\chi\right)\right\|_{H_h^{-N}(\Gamma)\to L^2(\Gamma)}$$
$$\lesssim x_1^{-\ell}h^{\ell}\theta^{-\ell-\widetilde k/2}\left(1+h^{\mu/2}\log(h^{-1})\theta^{-d-\mu}\right)e^{-Cx_1\theta/2h},
\end{equation}
where $\widetilde k=0$ if
 $k=-1$ and $\widetilde k=k$ if $k\ge 0$. 
If $\widetilde\chi\in S^p(\Gamma)$, ${\rm supp}\,\widetilde\chi\subset{\rm supp}(1-\eta)$, we have the bound
\begin{equation}\label{eq:4.14}
\left\|{\rm Op}_h\left(e^{ix_1\rho/h}\rho^{-k}\widetilde\chi\right)\right\|_{H_h^1(\Gamma)\to L^2(\Gamma)}\lesssim 
h^{\ell}x_1^{-\ell},
\end{equation}
provided $p<k+\ell+1$.
\end{prop}

{\it Proof.} We will first prove the following

\begin{lemma} Let $k$ and $\widetilde k$ be as in Proposition 4.4. Then, there exists a constant $C>0$ such that we have the bounds
\begin{equation}\label{eq:4.15} 
\left|\partial_{\xi'}^\alpha\partial_{x'}^\beta\left(e^{ix_1\rho/h}\rho^{-k}\right)\right|\le 
C_{\alpha}\theta^{-\widetilde k/2-|\alpha|-|\beta|}e^{-Cx_1\theta/h}
\quad\mbox{on}\quad{\rm supp}\,\eta,
\end{equation}
\begin{equation}\label{eq:4.16} 
\left|\partial_{\xi'}^\alpha\left(e^{ix_1\rho/h}\rho^{-k}\right)\right|\le 
C_{\alpha}|\xi'|^{-k-|\alpha|}e^{-Cx_1|\xi'|/h}
\quad\mbox{on}\quad{\rm supp}(1-\eta),
\end{equation}
for all multi-indices $\alpha$ and $\beta$ such that $|\beta|\le \mu$. 
\end{lemma}

{\it Proof.} It follows from Lemma 3.2 of \cite{kn:V1} that, if $k\ge 1$ is an integer or $k=-1$, we have the bounds 
\begin{equation}\label{eq:4.17} 
\left|\partial_{\xi'}^\alpha\partial_{x'}^\beta(\rho^{-k})\right|\le 
\left\{
\begin{array}{l}
C_{\alpha}|\rho|^{-k-2|\alpha|-2|\beta|}
\quad\mbox{on}\quad{\rm supp}\,\eta,\\
C_{\alpha}|\rho|^{-k-|\alpha|}
\quad\mbox{on}\quad{\rm supp}(1-\eta),
\end{array}
\right.
\end{equation}
for all multi-indices $\alpha$ and $\beta$ such that $|\beta|\le \mu$. By (\ref{eq:4.17}) together with (\ref{eq:4.2})
and (\ref{eq:4.3}) we get that the function $\rho^{-k}$, $k\ge 1$ being integer, satisfies the bounds
\begin{equation}\label{eq:4.18} 
\left|\partial_{\xi'}^\alpha\partial_{x'}^\beta(\rho^{-k})\right|\le 
\left\{
\begin{array}{l}
C_{\alpha}\theta^{-k/2-|\alpha|-|\beta|}
\quad\mbox{on}\quad{\rm supp}\,\eta,\\
C_{\alpha}|\xi'|^{-k-|\alpha|}
\quad\mbox{on}\quad{\rm supp}(1-\eta),
\end{array}
\right.
\end{equation}
while the function $\rho$ satisfies
\begin{equation}\label{eq:4.19} 
\left|\partial_{\xi'}^\alpha\partial_{x'}^\beta\rho\right|\le 
\left\{
\begin{array}{l}
C_{\alpha}|\rho|^{1-2|\alpha|-2|\beta|}
\quad\mbox{on}\quad{\rm supp}\,\eta,\\
C_{\alpha}|\xi'|^{1-|\alpha|}
\quad\mbox{on}\quad{\rm supp}(1-\eta),
\end{array}
\right.
\end{equation}
for all multi-indices $\alpha$ and $\beta$ such that $|\beta|\le \mu$. Clearly, the bounds (\ref{eq:4.18})
hold with $k=-1$, provided $|\alpha|+|\beta|\ge 1$. Therefore, to prove the lemma it suffices to prove the bounds 
(\ref{eq:4.15}) and (\ref{eq:4.16}) with $k=0$. This in turn follows from Lemma 4.2 of \cite{kn:V4}, but we will sketch the proof here
for the sake of completeness. 

Let us see that the functions 
$$c_{\alpha,\beta}=e^{-ix_1\rho/h}\partial_{\xi'}^\alpha
\partial_{x'}^\beta\left(e^{ix_1\rho/h}\right),\quad|\alpha|+|\beta|\ge 1,\,|\beta|\le\mu,$$
satisfy the bounds
\begin{equation}\label{eq:4.20} 
\left|\partial_{\xi'}^{\alpha'}\partial_{x'}^{\beta'}c_{\alpha,\beta}\right|\lesssim\sum_{j=1}^{|\alpha|+|\beta|+|\alpha'|+|\beta'|}
\left(\frac{x_1}{h|\rho|}\right)^j|\rho|^{-2(|\alpha|+|\beta|+|\alpha'|+|\beta'|-j)}
\end{equation}
on ${\rm supp}\,\eta$, and 
\begin{equation}\label{eq:4.21} 
\left|\partial_{\xi'}^{\alpha'}\partial_{x'}^{\beta'}c_{\alpha,\beta}\right|\lesssim\sum_{j=1}^{|\alpha|+|\beta|+|\alpha'|+|\beta'|}
\left(\frac{x_1}{h}\right)^j|\xi'|^{-(|\alpha|+|\alpha'|-j)}
\end{equation}
on ${\rm supp}(1-\eta)$, 
for all multi-indices $\alpha'$ and $\beta'$ such that $|\beta|+|\beta'|\le\mu$. We will proceed by induction in $|\alpha|+|\beta|$. 
Let $\alpha_1$ and $\beta_1$
be multi-indices such that $|\alpha_1|+|\beta_1|=1$ and observe that 
$$c_{\alpha+\alpha_1,\beta+\beta_1}=\partial_{\xi'}^{\alpha_1}\partial_{x'}^{\beta_1}c_{\alpha,\beta}
+ix_1h^{-1}c_{\alpha,\beta}\partial_{\xi'}^{\alpha_1}\partial_{x'}^{\beta_1}\rho.$$
More generally, we have 
\begin{equation}\label{eq:4.22} 
\partial_{\xi'}^{\alpha'}\partial_{x'}^{\beta'}c_{\alpha+\alpha_1,\beta+\beta_1}=
\partial_{\xi'}^{\alpha_1+\alpha'}\partial_{x'}^{\beta_1+\beta'}c_{\alpha,\beta}+ix_1h^{-1}\partial_{\xi'}^{\alpha'}\partial_{x'}^{\beta'}\left(c_{\alpha,\beta}\partial_{\xi'}^{\alpha_1}\partial_{x'}^{\beta_1}\rho\right).
\end{equation}
By (\ref{eq:4.19}) and (\ref{eq:4.22}), it is easy to see that if (\ref{eq:4.20})
and (\ref{eq:4.21}) hold for $c_{\alpha,\beta}$, they hold for $c_{\alpha+\alpha_1,\beta+\beta_1}$ as well. 

Using (\ref{eq:4.20}) together with (\ref{eq:4.2}) we obtain
$$\left|e^{ix_1\rho/h}c_{\alpha,\beta}\right|\lesssim
\sum_{j=1}^{|\alpha|+|\beta|}\left(\frac{x_1}{h|\rho|}\right)^j|\rho|^{-2(|\alpha|+|\beta|-j)}e^{-2C\theta x_1(h|\rho|)^{-1}}$$
 $$\lesssim\sum_{j=1}^{|\alpha|+|\beta|}\theta^{-j}|\rho|^{-2(|\alpha|+|\beta|-j)}e^{-Cx_1\theta/h}
 \lesssim\theta^{-|\alpha|-|\beta|}e^{-Cx_1\theta/h}.$$
Similarly, by (\ref{eq:4.21}) we obtain
$$\left|e^{ix_1\rho/h}c_{\alpha,\beta}\right|\lesssim\sum_{j=1}^{|\alpha|+|\beta|}
\left(\frac{x_1}{h}\right)^j|\xi'|^{-|\alpha|+j}e^{-2Cx_1|\xi'|/h}
\lesssim|\xi'|^{-|\alpha|}e^{-Cx_1|\xi'|/h}.$$
\eproof

In view of the inequality
$$e^{-Cx_1\theta/2h}\lesssim x_1^{-\ell}h^\ell\theta^{-\ell},$$
the bound (\ref{eq:4.13}) follows from Proposition 2.4 and (\ref{eq:4.15}). 
Since
$$ e^{-Cx_1|\xi'|/h}\lesssim h^\ell x_1^{-\ell}|\xi'|^{-\ell},$$
the bound (\ref{eq:4.14}) follows from Proposition 2.3 and (\ref{eq:4.16}).
\eproof

Since the function $A_1$ is supported in $\delta/2\le x_1\le \delta$, the next lemma is an immediate consequence of 
Lemma 4.5 with $k=-1$.

\begin{lemma} For every $m\ge 0$ we have the estimates
\begin{equation}\label{eq:4.23} 
\left|\partial_{\xi'}^\alpha\left(e^{ix_1\rho/h}A_1\right)\right|\le 
\left\{
\begin{array}{l}
C_{\alpha,m}h^m\theta^{-m-|\alpha|}
\quad\mbox{on}\quad{\rm supp}\,\eta,\\
C_{\alpha,m}h^m|\xi'|^{-m-|\alpha|}
\quad\mbox{on}\quad{\rm supp}(1-\eta),
\end{array}
\right.
\end{equation}
for all multi-indices $\alpha$ with constants $C_{\alpha,m}>0$ independent of $x_1$, $\theta$, $z$ and $h$.
\end{lemma}

This lemma together with Proposition 2.3 imply the following

\begin{prop}  Given any $N\ge 0$ independent of $h$, we have the estimate
\begin{equation}\label{eq:4.24}
\left\|\left({\rm Op}_h\left(e^{ix_1\rho/h}A_1\right)f\right)(x_1,\cdot)\right\|_{L^2(\Gamma)}\lesssim h\|f\|_{H_h^{-N}(\Gamma)}.
\end{equation}
\end{prop}

Let $u$ satisfy equation (\ref{eq:4.1}) with $u|_\Gamma={\rm Op}_h(\chi)f$. Then $u-\widetilde u$ satisfies equation (\ref{eq:3.1}) with
$v$ replaced by $\widetilde v$. Therefore, taking into account that 
$$-ih\partial_{\nu}\widetilde u|_{x_1=0}={\rm Op}_h(\rho\chi)f,$$
by (\ref{eq:3.2}) we get the estimate
\begin{equation}\label{eq:4.25}
\left\|{\cal N}(z,h){\rm Op}_h(\chi)f-{\rm Op}_h(\rho\chi)f\right\|_{L^2(\Gamma)}\lesssim h^{1/2}\theta^{-1/2}\|\widetilde v\|_{L^2(\Omega)}.
\end{equation}
On the other hand, Propositions 4.3 and 4.7 imply the following

\begin{lemma} If $\chi$ is of compact support, given any $N\ge 0$ independent of $h$, we have the estimate
\begin{equation}\label{eq:4.26}
\|\widetilde v\|_{L^2(\Omega)}\lesssim h^{1/2}\theta^{-2}\left(1+h^{\mu/2}\log(h^{-1})\theta^{-d-\mu}\right)\|f\|_{H_h^{-N}(\Gamma)}.
\end{equation}
If ${\rm supp}\,\chi\subset{\rm supp}(1-\eta)$, given any $0<\epsilon<1/2$ independent of $h$, we have the estimate
\begin{equation}\label{eq:4.27}
\|\widetilde v\|_{L^2(\Omega)}\lesssim h^{(1-3\epsilon)/2}\|f\|_{H_h^1(\Gamma)}.
\end{equation}
\end{lemma}

{\it Proof.} Clearly,
$$\|\widetilde v\|^2_{L^2(\Omega)}\lesssim \int_{\delta/2}^{\delta}
\left\|{\rm Op}_h\left(e^{ix_1\rho/h}A_1\right)f\right\|^2_{L^2(\Gamma)}dx_1+
\int_0^{\delta}\left\|{\rm Op}_h\left(e^{ix_1\rho/h}A_2\right)f\right\|^2_{L^2(\Gamma)}dx_1.$$
Therefore, if $\chi$ is of compact support, by (\ref{eq:4.11}) and (\ref{eq:4.24}), we get
$$\|\widetilde v\|^2_{L^2(\Omega)}\lesssim h\|f\|^2_{H_h^{-N}(\Gamma)}\int_{\delta/2}^{\delta}dx_1$$
$$+\theta^{-3}\left(1+h^{\mu/2}\log(h^{-1})\theta^{-d-\mu}\right)^2\|f\|^2_{H_h^{-N}(\Gamma)}
\int_0^{\delta}e^{-Cx_1\theta/h}dx_1$$
$$\lesssim h\theta^{-4}\left(1+h^{\mu/2}\log(h^{-1})\theta^{-d-\mu}\right)^2\|f\|^2_{H_h^{-N}(\Gamma)}$$
for every $N\ge 0$, which clearly implies (\ref{eq:4.26}). Similarly, if ${\rm supp}\,\chi\subset{\rm supp}(1-\eta)$, 
by (\ref{eq:4.12}) and (\ref{eq:4.24}), we get
$$\|\widetilde v\|^2_{L^2(\Omega)}\lesssim h\|f\|^2_{H_h^{-N}(\Gamma)}\int_{\delta/2}^{\delta}dx_1+h^{1-3\epsilon}\|f\|^2_{H_h^1(\Gamma)}
\int_0^{\delta}x_1^{-1+2\epsilon}dx_1\lesssim h^{1-3\epsilon}\|f\|^2_{H_h^1(\Gamma)}$$
which implies (\ref{eq:4.27}).
\eproof

Combining (\ref{eq:4.25}) with Lemma 4.8 (with $\epsilon$ small enough) and taking into account that $\theta^{-1/2}\le h^{-1/5}$,
we obtain the following

\begin{lemma} If $\chi$ is of compact support, given any $N\ge 0$ independent of $h$, we have the estimate
\begin{equation}\label{eq:4.28}
\left\|{\cal N}(z,h){\rm Op}_h(\chi)f-{\rm Op}_h(\rho\chi)f\right\|_{L^2(\Gamma)}\lesssim h\theta^{-5/2}
\left(1+h^{\mu/2}\log(h^{-1})\theta^{-d-\mu}\right)\|f\|_{H_h^{-N}(\Gamma)}.
\end{equation}
If ${\rm supp}\,\chi\subset{\rm supp}(1-\eta)$, we have the estimate
\begin{equation}\label{eq:4.29}
\left\|{\cal N}(z,h){\rm Op}_h(\chi)f-{\rm Op}_h(\rho\chi)f\right\|_{L^2(\Gamma)}\lesssim h^{3/5}\|f\|_{H_h^1(\Gamma)}.
\end{equation}
\end{lemma}

Since the function $\eta$ can be written as a finite sum of functions $\chi$ for which (\ref{eq:4.28}) holds and 
the function $1-\eta$ can be written as a finite sum of functions $\chi$ for which (\ref{eq:4.29}) holds, Lemma 4.9 implies the following

\begin{prop} Given any $N\ge 0$ independent of $h$, we have the estimates
\begin{equation}\label{eq:4.30}
\left\|{\cal N}(z,h){\rm Op}_h(\eta)f-{\rm Op}_h(\rho\eta)f\right\|_{L^2(\Gamma)}\lesssim h\theta^{-5/2}
\left(1+h^{\mu/2}\log(h^{-1})\theta^{-d-\mu}\right)\|f\|_{H_h^{-N}(\Gamma)}.
\end{equation}
and
\begin{equation}\label{eq:4.31}
\left\|{\cal N}(z,h){\rm Op}_h(1-\eta)f-{\rm Op}_h(\rho(1-\eta))f\right\|_{L^2(\Gamma)}\lesssim h^{3/5}\|f\|_{H_h^1(\Gamma)}.
\end{equation}
\end{prop}

Clearly, the estimate (\ref{eq:4.4}) follows from (\ref{eq:4.30}) and (\ref{eq:4.31}). To prove the estimate (\ref{eq:4.5}) we need
to improve (\ref{eq:4.31}), only, assuming that $c\equiv 1$. To this end, we have to build a better parametrix
in the elliptic region in this case. This will be carried out in the next section. 

\section{Improved parametrix in the elliptic region}

Our goal in this section is to prove the following

\begin{prop} Suppose that $c\equiv 1$ and $n\in C^1(\Omega_\delta)$, $n|_\Gamma\in C^2(\Gamma)$. Then we have the estimate
\begin{equation}\label{eq:5.1}
\left\|{\cal N}(z,h){\rm Op}_h(1-\eta)f-{\rm Op}_h(\rho(1-\eta)+hq)f\right\|_{L^2(\Gamma)}\lesssim h\|f\|_{H_h^{-1}(\Gamma)},
\end{equation}
where $q\in S^0(\Gamma)$ is independent of $z$, $h$ and the function $n$.
\end{prop}

Then the estimate (\ref{eq:4.5}) would follow from (\ref{eq:4.30}) and (\ref{eq:5.1}). To prove (\ref{eq:5.1}) we will improve our
parametrix when  $\chi\in S^0(\Gamma)$ with ${\rm supp}\,\chi\subset{\rm supp}(1-\eta)$. To this end, we choose the phase function in the form
$$\varphi=-\langle x',\xi'\rangle+x_1\rho+x_1^2\varphi_2+x_1^3\varphi_3,$$
where the functions $\varphi_2,\varphi_3\in S^1(\Gamma)$ do not depend on $x_1$, $z$ and $h$, and will be choosen in such a way that the eikonal equation (\ref{eq:4.6})
is satisfied with $\Phi={\cal O}\left(1+x_1^2|\xi'|^2\right)$. Furthermore, we choose the amplitude in the form
$$a=a_0+x_1a_1+x_1^2a_2,$$
where $a_0=\chi$, $a_1=a_{1,0}+ha_{1,1}$ with functions $a_{1,0}, a_2\in S^0(\Gamma)$, $a_{1,1}\in S^{-1}(\Gamma)$ 
independent of $x_1$ and $h$, and will be choosen in such a way that
$A_3={\cal O}\left(x_1^2|\xi'|+hx_1\right)$. 
To find $\varphi_2$ and $\varphi_3$, observe that the function $\Phi$ 
 with the new phase can be written in the form
$$\Phi=2\rho(2\varphi_2+3x_1\varphi_3)+x_1(2\varphi_2+3x_1\varphi_3)^2-zn^\sharp$$
$$-2\langle R_0(x')\xi',\nabla_{x'}(\rho+x_1\varphi_2+x_1^2\varphi_3)\rangle$$
$$+x_1\langle R_0(x')\nabla_{x'}(\rho+x_1\varphi_2+x_1^2\varphi_3),\nabla_{x'}(\rho+x_1\varphi_2+x_1^2\varphi_3)\rangle$$
$$+\langle R^\sharp(x)\xi',\xi'\rangle-2x_1\langle R^\sharp(x)\xi',\nabla_{x'}(\rho+x_1\varphi_2+x_1^2\varphi_3)\rangle$$
$$+x_1^2\langle R^\sharp(x)\nabla_{x'}(\rho+x_1\varphi_2+x_1^2\varphi_3),\nabla_{x'}(\rho+x_1\varphi_2+x_1^2\varphi_3)\rangle.$$
We need now the following

\begin{lemma} On ${\rm supp}(1-\eta)$ we have the bounds
\begin{equation}\label{eq:5.2} 
\left|\partial_{\xi'}^\alpha\partial_{x'}^\beta\left(\rho-ir_0^{1/2}\right)\right|\le 
C_{\alpha}|\xi'|^{-1-|\alpha|},
\end{equation} 
\begin{equation}\label{eq:5.3} 
\left|\partial_{\xi'}^\alpha\left(\rho^{-1}+ir_0^{-1/2}\right)\right|\le 
C_{\alpha}|\xi'|^{-3-|\alpha|},
\end{equation}
\begin{equation}\label{eq:5.4} 
\left|\partial_{\xi'}^\alpha\left(\nabla_{x'}\rho-i2^{-1}r_0^{-1/2}\nabla_{x'}r_0\right)\right|\le 
C_{\alpha}|\xi'|^{-1-|\alpha|},
\end{equation}
 for all multi-indices $\alpha$ and $\beta$ such that $|\beta|\le 2$.
 \end{lemma}

{\it Proof.} The bound (\ref{eq:5.2}) with $\beta=0$ follows from the identity
$$\rho-ir_0^{1/2}=\frac{zn_0}{\rho+ir_0^{1/2}}$$
together with (\ref{eq:4.19}) and the fact that $|\rho+ir_0^{1/2}|\ge Cr_0^{1/2}$, $C>0$, on ${\rm supp}(1-\eta)$. 
Since $n_0\in C^2$, we can differentiate the above identity with respect to $x'$ to get (\ref{eq:5.2}) with $1\le|\beta|\le 2$
in the same way. Furthermore, since 
$$\rho^{-1}+ir_0^{-1/2}=i\rho^{-1}r_0^{-1/2}\left(\rho-ir_0^{1/2}\right),$$
the bound (\ref{eq:5.3}) follows from (\ref{eq:5.2}) and (\ref{eq:4.18}). Finally,
the bound (\ref{eq:5.4}) follows from the identity (\ref{eq:4.7}) together with the bounds (\ref{eq:5.3}) and (\ref{eq:4.18}).
\eproof

Write the function $R^\sharp$ as
$$R^\sharp(x)=R^\sharp_0(x')+x_1R^\sharp_1(x')+x_1^2R^\flat(x),$$
where $R^\sharp_0$, $R^\sharp_1$, $R^\flat$ are smooth functions. 
Then the function $\Phi$ can be written in the form
$$\Phi=\Phi_0+x_1\Phi_1+{\cal O}\left(1+x_1^2|\xi'|^2\right),$$
where 
$$\Phi_0=4ir_0^{1/2}\varphi_2-ir_0^{-1/2}\langle R_0(x')\xi',\nabla_{x'}r_0\rangle+\langle R^\sharp_0(x')\xi',\xi'\rangle,$$
$$\Phi_1=6ir_0^{1/2}\varphi_3+4\varphi_2^2-2\langle R_0(x')\xi',\nabla_{x'}\varphi_2\rangle+\langle R^\sharp_1(x')\xi',\xi'\rangle$$
$$-(4r_0)^{-1}\langle R_0(x')\nabla_{x'}r_0,\nabla_{x'}r_0\rangle-ir_0^{-1/2}\langle R^\sharp_0(x')\xi',\nabla_{x'}r_0\rangle.$$
We now choose the function $\varphi_2$ so that $\Phi_0=0$ and the function $\varphi_3$ so that $\Phi_1=0$.
Clearly, the functions $\varphi_2$ and $\varphi_3$ are smooth and satisfy the estimates
\begin{equation}\label{eq:5.5} 
\left|\partial_{\xi'}^\alpha\partial_{x'}^\beta\varphi_j\right|\le 
C_{\alpha,\beta}|\xi'|^{1-|\alpha|},\quad j=2,3,
\end{equation}
on supp$\chi$, for all multi-indices $\alpha$ and $\beta$.  It is easy to see that (\ref{eq:5.5}) together with
Lemma 5.2 imply the following

\begin{lemma} For $0<x_1\le \delta$, we have the bounds
\begin{equation}\label{eq:5.6} 
\left|\partial_{\xi'}^\alpha\Phi\right|\le 
C_{\alpha}\left(1+x_1^2|\xi'|^2\right)|\xi'|^{-|\alpha|}
\end{equation}
for all multi-indices $\alpha$ with constants $C_{\alpha}>0$ independent of $x_1$, $z$ and $h$. 
\end{lemma}

To find the functions $a_1$ and $a_2$, observe that in this case the function $A_3$ is of the form
$$A_3=-2i(\rho+2x_1\varphi_2+3x_1^2\varphi_3)(a_1+2x_1a_2)$$
$$-2i\langle-\xi'+x_1\nabla_{x'}(\rho+x_1\varphi_2+x_1^2\varphi_3),\nabla_{x'}(\chi+x_1a_1+x_1^2a_2)\rangle$$
$$-2i(\chi+x_1a_1+x_1^2a_2)(\varphi_2+x_1\varphi_3)$$
$$-ix_1(\chi+x_1a_1+x_1^2a_2)\left\langle R(x)\nabla_{x'},\nabla_{x'}\right\rangle(\rho+x_1\varphi_2+x_1^2\varphi_3)$$
$$-h\left\langle R(x)\nabla_{x'},\nabla_{x'}\right\rangle(\chi+x_1a_1+x_1^2a_2)$$
$$+Q_1(\chi+x_1a_1+x_1^2a_2)(\rho+2x_1\varphi_2+3x_1^2\varphi_3)$$
$$+\left\langle\widetilde Q,-\xi'+x_1\nabla_{x'}
(\rho+x_1\varphi_2)\right\rangle(\chi+x_1a_1+x_1^2a_2)$$ 
$$-ih\left\langle\widetilde Q,\nabla_{x'}(\chi+x_1a_1+x_1^2a_2)\right\rangle.$$
Since $c\equiv 1$, the functions $Q_1$ and $\widetilde Q$ are smooth, so we can write them in the form
$$Q_1(x)=Q_{1,0}(x')+x_1Q_{1,1}(x')+x_1^2Q^\flat_1(x),$$
$$\widetilde Q(x)=\widetilde Q_0(x')+x_1\widetilde Q_1(x')+x_1^2\widetilde Q^\flat(x),$$
where all functions are smooth. Hence we can write the function $A_3$ in the form
$$A_3=A_{3,0}+x_1A_{3,1}+{\cal O}\left(x_1^2|\xi'|+hx_1\right),$$
where
$$A_{3,0}=2r_0^{1/2}a_1+2i\langle\xi',\nabla_{x'}\chi\rangle-2i\chi\varphi_2+iQ_1\chi r_0^{1/2}
-\left\langle\widetilde Q,\xi'\right\rangle\chi$$
$$-h\left\langle R_0(x')\nabla_{x'},\nabla_{x'}\right\rangle\chi-ih\left\langle\widetilde Q,\nabla_{x'}\chi\right\rangle,$$ 
$$A_{3,1}=4r_0^{1/2}a_2-4ia_1\varphi_2+2i\langle\xi',\nabla_{x'}a_{1,0}\rangle
+r_0^{-1/2}\langle\nabla_{x'}r_0,\nabla_{x'}\chi\rangle$$
$$-2i\chi\varphi_3-2ia_1\varphi_2+\chi\left\langle R_0(x')\nabla_{x'},\nabla_{x'}\right\rangle r_0^{1/2}$$
$$+2\chi Q_{1,0}(x')\varphi_2+iQ_{1,0}(x')r_0^{1/2}a_{1,0}+i\chi Q_{1,1}(x')r_0^{1/2}$$
$$-\chi\left\langle\widetilde Q_1(x'),\xi'\right\rangle
-\left\langle\widetilde Q_0(x'),\xi'\right\rangle a_{1,0}
+i2^{-1}r_0^{-1/2}\left\langle\widetilde Q_0(x'),\nabla_{x'}r_0\right\rangle\chi.$$ 

We now choose the function $a_1$ so that $A_{3,0}=0$ and the function $a_2$ so that $A_{3,1}=0$. We get
$$a_{1,0}=-ir_0^{-1/2}\langle\xi',\nabla_{x'}\chi\rangle+ir_0^{-1/2}\chi\varphi_2-i2^{-1}Q_1\chi 
+2^{-1}r_0^{-1/2}\left\langle\widetilde Q,\xi'\right\rangle\chi\in S^0(\Gamma),$$
$$a_{1,1}=2^{-1}r_0^{-1/2}\left\langle R_0(x')\nabla_{x'},\nabla_{x'}\right\rangle\chi+i2^{-1}r_0^{-1/2}\left\langle\widetilde Q,\nabla_{x'}\chi\right\rangle\in S^{-1}(\Gamma).$$ 
 The next lemma follows easily from (\ref{eq:5.5}) and Lemma 5.2.
 
\begin{lemma} For $0<x_1\le \delta$, we have the bounds
\begin{equation}\label{eq:5.7} 
\left|\partial_{\xi'}^\alpha A_3\right|\le 
C_{\alpha}\left(hx_1+x_1^2|\xi'|\right)|\xi'|^{-|\alpha|}
\end{equation}
for all multi-indices $\alpha$ with constants $C_{\alpha}>0$ independent of $x_1$, $z$ and $h$. 
\end{lemma}
 
Set 
$$\widetilde\rho=\rho+x_1\varphi_2+x_1^2\varphi_3=\rho+{\cal O}(x_1|\xi'|).$$
Taking $\delta$ small enough we can arrange that the inequalities in (\ref{eq:4.3}) still hold with $\rho$
replaced by $\widetilde\rho$ for all $0<x_1\le\delta$. Therefore, the bound (\ref{eq:4.16}) (with $k=0$) holds with 
$e^{ix_1\rho/h}$ replaced by $e^{ix_1\widetilde\rho/h}$. This together with Lemmas 5.3 and 5.4 imply the following 

\begin{lemma} For $0<x_1\le \delta$ and for all $0<\epsilon<1/2$, we have the estimates
\begin{equation}\label{eq:5.8} 
\left|\partial_{\xi'}^\alpha\left(e^{ix_1\widetilde\rho/h}A_2\right)\right|\le 
C_{\alpha}h^{5/2-\epsilon}x_1^{-1/2+\epsilon}|\xi'|^{-3/2+\epsilon-|\alpha|}
\end{equation}
for all multi-indices $\alpha$ with constants $C_{\alpha}>0$ independent of $x_1$, $z$ and $h$. 
\end{lemma}

It is easy also to see that the new function $A_1$ satisfies the bound (\ref{eq:4.23}) (on supp$(1-\eta)$) 
with 
$e^{ix_1\rho/h}$ replaced by $e^{ix_1\widetilde\rho/h}$. 
Hence Proposition 4.7 still holds with 
$e^{ix_1\rho/h}$ replaced by $e^{ix_1\widetilde\rho/h}$. 
Therefore, in the same way as in the previous section we can deduce from Lemma 5.5 and Proposition 2.3 that the function $\widetilde v$ satisfies
the estimate
\begin{equation}\label{eq:5.9}
\|\widetilde v\|_{L^2(\Omega)}\lesssim h^{(5-3\epsilon)/2}\|f\|_{H_h^{-1}(\Gamma)}.
\end{equation}
 We also have
$$-ih\partial_{\nu}\widetilde u|_{x_1=0}={\rm Op}_h(\rho\chi-iha_1)f={\rm Op}_h(\rho\chi-iha_{1,0})f-ih^2{\rm Op}_h(a_{1,1})f$$
and the operator ${\rm Op}_h(a_{1,1}):H_h^{-1}(\Gamma)\to L^2(\Gamma)$ is uniformly bounded. 
Therefore, Proposition 5.1 follows from (\ref{eq:3.2}) and  (\ref{eq:5.9}). 

\section{Proof of Theorem 1.1}

We have to show that if $\lambda\in \mathbb{C}$ belongs to the eigenvalue-free regions in Theorem 1.1, then the solution
$(u_1,u_2)$ of the equation (\ref{eq:1.1}) is identically zero. Clearly, it suffices to show that the function
$f=u_1|_\Gamma=u_2|_\Gamma$ is identically zero. We may suppose that $|\lambda|\gg 1$. Set $h=|\lambda|^{-1}\ll 1$ 
and $z=(h\lambda)^2$. Clearly, $z\in Z^+$ if ${\rm Re}\,\lambda\ge |{\rm Im}\,\lambda|$ and $z\in Z^-$ if 
$|{\rm Im}\,\lambda|\ge {\rm Re}\,\lambda$, where $Z^+$ and $Z^-$ are as in Section 3.  
Set $\theta=|{\rm Im}\,z|$ if $z\in Z^+$ and $\theta=1$ if $z\in Z^-$.  
In order to simplify the notations, in what follows the restrictions of the functions 
$c_j,n_j,\widetilde n_j$, $j=1,2$, on $\Gamma$ will be again denoted by $c_j,n_j,\widetilde n_j$, respectively. 
By assumption, we have that they belong to $C^\mu(\Gamma)$. Define the functions $\rho_j$, $j=1,2$, by replacing in the definition of $\rho$
the function $\widetilde n$ by $\widetilde n_j$. 
Denote by ${\cal N}_j$, $j=1,2$, the Dirichlet-to-Neumann map 
introduced in Section 4 associated to $(c_j,n_j)$, and set
$$T(z,h)=c_1{\cal N}_1(z,h)-c_2{\cal N}_2(z,h).$$
If $\lambda$ is a transmission eigenvalue, then $Tf=0$ on $\Gamma$. By Theorem 4.1 we get
\begin{equation}\label{eq:6.1}
\left\|{\rm Op}_h(c_1\rho_1-c_2\rho_2)f\right\|_{L^2(\Gamma)}\lesssim h^{3/5}\|f\|_{H_h^k(\Gamma)}$$  $$+h\theta^{-5/2}\left(1+h^{\mu/2}\log(h^{-1})\theta^{-d-\mu}\right)\|f\|_{H_h^k(\Gamma)}
\end{equation}
provided $\theta\ge h^{2/5}$, where $k=1$ in the anisotropic case and $k=-1$ in the isotropic case.
Set
$$a_1=(r_0+1)^{-1/2}(c_1\rho_1+c_2\rho_2),\quad a_2=c_1\rho_1-c_2\rho_2.$$
We will prove now  the following

\begin{prop} Suppose that $\theta$ satisfes the condition 
\begin{equation}\label{eq:6.2}
h^{\mu/2}\log(h^{-1})\theta^{1/2-d-\mu}\le 1.
\end{equation}
Then we have the bounds
\begin{equation}\label{eq:6.3}
\left\|{\rm Op}_h(a_1)\right\|_{L^2(\Gamma)\to L^2(\Gamma)}
 \lesssim 1, 
\end{equation}
\begin{equation}\label{eq:6.4}
\left\|{\rm Op}_h(a_1){\rm Op}_h(a_2)-{\rm Op}_h(a_1a_2)\right\|_{H_h^{k}(\Gamma)\to L^2(\Gamma)}
 \lesssim h\theta^{-1}+h^{\mu/2}\log(h^{-1})\theta^{1/2-d-\mu}. 
\end{equation}
\end{prop} 

{\it Proof.} Clearly, the functions $a_1$ and $a_2$ are bounded on supp$\,\eta$. Moreover, 
it follows from (\ref{eq:4.19}) and (\ref{eq:4.2}) that they satisfy the bounds
\begin{equation}\label{eq:6.5} 
\left|\partial_{\xi'}^\alpha\partial_{x'}^\beta a_j\right|\lesssim\theta^{1/2-|\alpha|-|\beta|}
\quad\mbox{on}\quad{\rm supp}\,\eta,
\end{equation}
for all multi-indices $\alpha$ and $\beta$ such that $|\alpha|+|\beta|\ge 1$ and $|\beta|\le \mu$.
On the other hand, since $\rho_j=ir_0^{1/2}+{\cal O}(r_0^{-1/2})$ for $r_0\gg 1$, we have
$$a_1=i(c_1+c_2)+{\cal O}(r_0^{-1}),$$
$$a_2=i(c_1-c_2)r_0^{1/2}+{\cal O}(r_0^{-1/2}),\quad\quad if\quad k=1,$$
$$a_2=\frac{z(n_1-n_2)}{\rho_1+\rho_2}=
(2i)^{-1}z(n_1-n_2)r_0^{-1/2}+{\cal O}(r_0^{-3/2}),\quad\quad if\quad k=-1.$$
Set $b_1=i(c_1+c_2)$ and 
$$b_2=\left\{
\begin{array}{l}
i(c_1-c_2)
\quad\mbox{if}\quad k=1,\\
(2i)^{-1}z(n_1-n_2)
\quad\mbox{if}\quad k=-1.
\end{array}
\right.
$$
We now decompose the functions $a_j$ as follows
$$a_j=\sum_{\ell=1}^3a_j^{(\ell)},$$
where $a_j^{(1)}=\eta a_j$, $a_1^{(3)}=b_1(1-\eta)$ and 
$a_2^{(3)}=b_2(1-\eta)r_0^{k/2}$. Then the functions $a_j^{(2)}$ are supported on 
supp$(1-\eta)$ and $a_1^{(2)}={\cal O}(r_0^{-1})$, $a_2^{(2)}={\cal O}(r_0^{k/2-1})$. Moreover, they satisfy the bounds
\begin{equation}\label{eq:6.6} 
\left|\partial_{\xi'}^\alpha\partial_{x'}^\beta a_1^{(2)}\right|\lesssim\langle\xi'\rangle^{-2-|\alpha|},
\end{equation}
\begin{equation}\label{eq:6.7} 
\left|\partial_{\xi'}^\alpha\partial_{x'}^\beta a_2^{(2)}\right|\lesssim\langle\xi'\rangle^{k-2-|\alpha|},
\end{equation}
for all multi-indices $\alpha$ and $\beta$ such that $|\beta|\le \mu$. These bounds follow from the fact that 
(\ref{eq:5.2}) holds for all $\beta$ such that $|\beta|\le \mu$, which in turn can be easily proved by induction in $\mu$. 

We can write
$${\rm Op}_h(a_1){\rm Op}_h(a_2)-{\rm Op}_h(a_1a_2)=\sum_{\ell_1=1}^3\sum_{\ell_2=1}^3
\left({\rm Op}_h(a_1^{(\ell_1)}){\rm Op}_h(a_2^{(\ell_2)})-{\rm Op}_h(a_1^{(\ell_1)}a_2^{(\ell_2)})\right)$$
$$=\sum_{\ell_1=1}^2\sum_{\ell_2=1}^2
\left({\rm Op}_h(a_1^{(\ell_1)}){\rm Op}_h(a_2^{(\ell_2)})-{\rm Op}_h(a_1^{(\ell_1)}a_2^{(\ell_2)})\right)$$
 $$+b_1\sum_{\ell_2=1}^2
\left({\rm Op}_h((1-\eta)){\rm Op}_h(a_2^{(\ell_2)})-{\rm Op}_h((1-\eta)a_2^{(\ell_2)})\right)$$
$$+b_2\sum_{\ell_1=1}^2
\left({\rm Op}_h(a_1^{(\ell_1)}){\rm Op}_h((1-\eta)r_0^{k/2})-{\rm Op}_h(a_1^{(\ell_1)}(1-\eta)r_0^{k/2})\right)$$
$$+\sum_{\ell_1=1}^2
\left[{\rm Op}_h(a_1^{(\ell_1)}),b_2\right]{\rm Op}_h((1-\eta)r_0^{k/2})$$
$$+b_1\left[{\rm Op}_h(1-\eta),b_2\right]{\rm Op}_h((1-\eta)r_0^{k/2})$$
$$+b_1b_2\left({\rm Op}_h(1-\eta){\rm Op}_h((1-\eta)r_0^{k/2})-{\rm Op}_h((1-\eta)^2r_0^{k/2})\right).$$
Now we apply Proposition 2.6 to the operators in the double sum in the right-hand side, Proposition 2.5 to the operators
in the next two sums and Proposition 2.2 to the operator in the last term. Furthermore, the operator 
$${\rm Op}_h((1-\eta)r_0^{k/2}):H_h^k(\Gamma)\to L^2(\Gamma)$$
is bounded, while the norm of the commutators 
$$\left[{\rm Op}_h(a_1^{(\ell_1)}),b_2\right]:L^2(\Gamma)\to L^2(\Gamma)$$
can be bounded by using (\ref{eq:2.26}), and  
the norm of the commutator 
$$\left[{\rm Op}_h(1-\eta),b_2\right]=-\left[{\rm Op}_h(\eta),b_2\right]:L^2(\Gamma)\to L^2(\Gamma)$$
can be bounded by using (\ref{eq:2.22}). 
Thus, taking into account the bounds
(\ref{eq:6.5}), (\ref{eq:6.6}) and (\ref{eq:6.7}), we arrive at (\ref{eq:6.4}). The bound (\ref{eq:6.3}) follows
similarly from Propositions 2.1 and 2.4.
\eproof

Combining Proposition 6.1 with the estimate (\ref{eq:6.1}) we get
\begin{equation}\label{eq:6.8}
\left\|{\rm Op}_h(a_1a_2)f\right\|_{L^2(\Gamma)}\lesssim h^{3/5}\|f\|_{H_h^k(\Gamma)}+{\cal E}_1(h,\theta)\|f\|_{H_h^k(\Gamma)}
\end{equation}
provided $\theta\ge h^{2/5}$ and $\theta$ satisfies (\ref{eq:6.2}), where
$${\cal E}_1=h\theta^{-5/2}\left(1+h^{\mu/2}\log(h^{-1})\theta^{-d-\mu}\right)
+h^{\mu/2}\log(h^{-1})\theta^{1/2-d-\mu}.$$
Set
$$\theta_1(h)=
\left\{
\begin{array}{l}
\left(h^{\mu/2}\log(h^{-1})\right)^{1/(d+\mu-1/2)}
\quad\mbox{if}\quad\mu\le 2d-1,\\
\left(h^{\mu/2+1}\log(h^{-1})\right)^{1/(d+\mu+5/2)}
\quad\mbox{if}\quad 2d-1<\mu\le 4d,\\
h^{2/5}
\quad\mbox{if}\quad \mu>4d. 
\end{array}
\right.
$$
It is easy to check that $\theta\gg\theta_1$ implies ${\cal E}_1\ll 1$ and $\theta\ge h^{2/5}$ satisfies (\ref{eq:6.2}). 
We would like to show that for such $\theta$
and $h$ small enough the estimate (\ref{eq:6.8}) implies $f\equiv 0$ in the isotropic case. To do so, we will use the identity
\begin{equation}\label{eq:6.9}
(c_1\rho_1+c_2\rho_2)(c_1\rho_1-c_2\rho_2)=c_1^2\rho_1^2-c_2^2\rho_2^2=-(c_1^2-c_2^2)r_0+z(c_1n_1-c_2n_2).
\end{equation}
So, in the isotropic case (\ref{eq:6.9}) gives
$${\rm Op}_h(a_1a_2)=z(n_1-n_2){\rm Op}_h\left((r_0+1)^{-1/2}\right).$$
Recall that in this case we have $|n_1-n_2|\ge c$ with some constant $c>0$ by assumption. 
Since 
$$\left\|{\rm Op}_h\left((r_0+1)^{-1/2}\right)f\right\|_{L^2(\Gamma)}\sim\|f\|_{H_h^{-1}(\Gamma)},$$
the estimate (\ref{eq:6.8}) implies 
\begin{equation}\label{eq:6.10}
\|f\|_{H_h^{-1}(\Gamma)}\lesssim h^{3/5}\|f\|_{H_h^{-1}(\Gamma)}+{\cal E}_1\|f\|_{H_h^{-1}(\Gamma)}.
\end{equation}
Taking $h$ and ${\cal E}_1$ small enough we deduce from (\ref{eq:6.10}) that $\|f\|_{H_h^{-1}(\Gamma)}=0$, and hence 
$f\equiv 0$, as desired. Hence the region $\theta\gg\theta_1$ is an eigenvalue-free region. It is easy to see that this region corresponds
to (\ref{eq:1.12}) on the $\lambda$ plane. 

Consider now the anisotropic case. Then the function
$$m=\frac{c_1n_1-c_2n_2}{c_1^2-c_2^2}$$
is strictly negative under the condition (\ref{eq:1.6}) and strictly positive under the condition (\ref{eq:1.8}).
Moreover, we have $m\in C^\mu(\Gamma)$. Set
$$A_1=(r_0+1)(r_0-zm)^{-1},\quad A_2=(r_0+1)^{-1/2}(r_0-zm).$$
Let $\tau=|{\rm Im}\,z|$ if $z\in Z^-$ and $\tau=1$ if $z\in Z^+$. 
We need now the following

\begin{prop} Suppose the condition (\ref{eq:1.8}) fulfilled. Suppose also that $\theta$ satisfes the condition 
\begin{equation}\label{eq:6.11}
h^{\mu/2}\log(h^{-1})\theta^{-d-\mu}\le 1.
\end{equation}
Then we have the bounds
\begin{equation}\label{eq:6.12}
\left\|{\rm Op}_h(A_1)\right\|_{L^2(\Gamma)\to L^2(\Gamma)}
 \lesssim \theta^{-1}, 
\end{equation}
\begin{equation}\label{eq:6.13}
\left\|{\rm Op}_h(A_1){\rm Op}_h(A_2)-{\rm Op}_h(A_1A_2)\right\|_{H_h^1(\Gamma)\to L^2(\Gamma)}
 \lesssim h\theta^{-2}+h^{\mu/2}\log(h^{-1})\theta^{-d-\mu-1}. 
\end{equation}
Suppose the condition (\ref{eq:1.6}) fulfilled. Suppose also that $\tau\ge h^{1/2}$ satisfes the condition 
\begin{equation}\label{eq:6.14}
h^{\mu/2}\log(h^{-1})\tau^{-d-\mu}\le 1.
\end{equation}
Then we have the bounds
\begin{equation}\label{eq:6.15}
\left\|{\rm Op}_h(A_1)\right\|_{L^2(\Gamma)\to L^2(\Gamma)}
 \lesssim \tau^{-1}, 
\end{equation}
\begin{equation}\label{eq:6.16}
\left\|{\rm Op}_h(A_1){\rm Op}_h(A_2)-{\rm Op}_h(A_1A_2)\right\|_{H_h^1(\Gamma)\to L^2(\Gamma)}
 \lesssim h\tau^{-2}+h^{\mu/2}\log(h^{-1})\tau^{-d-\mu-1}. 
\end{equation}
\end{prop} 

{\it Proof.} Clearly, the function $A_2$ satisfies the bounds
\begin{equation}\label{eq:6.17} 
\left|\partial_{\xi'}^\alpha\partial_{x'}^\beta A_2\right|\lesssim\langle\xi'\rangle^{1-|\alpha|},
\end{equation}
for all multi-indices $\alpha$ and $\beta$ such that $|\beta|\le \mu$. We would like to find similar bounds for
the derivatives of $A_1$. To this end, observe that we can arrange $|r_0-zm|\ge Cr_0$, $C>0$,
on supp$(1-\eta)$, provided the constant $C_0$ in the definition of $\eta$ is taken large enough
(what we can do without loss of generality). 
The function $A_1$ satisfies the following

\begin{lemma} If $m{\rm Re}\,z\ge 0$ we have the bounds
\begin{equation}\label{eq:6.18} 
\left|\partial_{\xi'}^\alpha\partial_{x'}^\beta A_1\right|\lesssim
\left\{
\begin{array}{l}
|{\rm Im}\,z|^{-1-|\alpha|-|\beta|}
\quad\mbox{on}\quad{\rm supp}\,\eta,\\
|\xi'|^{-|\alpha|}
\quad\mbox{on}\quad{\rm supp}(1-\eta),
\end{array}
\right.
\end{equation}
for all multi-indices $\alpha$ and $\beta$ such that $|\beta|\le \mu$.
 If $m{\rm Re}\,z\le 0$ the bounds (\ref{eq:6.18}) hold with $|{\rm Im}\,z|$ replaced by $1$. 
\end{lemma}

{\it Proof.} We have 
$$|r_0-zm|\ge C(r_0+1),$$
if $m{\rm Re}\,z\le 0$, and
$$|r_0-zm|\ge C|{\rm Im}\,z|\quad\mbox{on}\quad{\rm supp}\,\eta,$$
$$|r_0-zm|\ge C|\xi'|^2\quad\mbox{on}\quad{\rm supp}(1-\eta),$$
if $m{\rm Re}\,z\ge 0$, where $C>0$ is a constant. The bounds (\ref{eq:6.18}) can be easily derived from the above
inequalities by induction in $|\alpha|+|\beta|$. 
\eproof

Decompose the functions $A_j$ as follows
$$A_j=\sum_{\ell=1}^3A_j^{(\ell)},$$
where $A_j^{(1)}=\eta A_j$, $A_1^{(3)}=(1-\eta)$ and 
$A_2^{(3)}=(1-\eta)r_0^{1/2}$. Then the functions $A_j^{(2)}$ are supported on 
supp$(1-\eta)$ and $A_1^{(2)}={\cal O}(r_0^{-1})$, $A_2^{(2)}={\cal O}(r_0^{-1/2})$. Moreover, it follows from
Lemma 6.3 that they satisfy the bounds
\begin{equation}\label{eq:6.19} 
\left|\partial_{\xi'}^\alpha\partial_{x'}^\beta A_1^{(2)}\right|\lesssim\langle\xi'\rangle^{-2-|\alpha|},
\end{equation}
\begin{equation}\label{eq:6.20} 
\left|\partial_{\xi'}^\alpha\partial_{x'}^\beta A_2^{(2)}\right|\lesssim\langle\xi'\rangle^{-1-|\alpha|},
\end{equation}
for all multi-indices $\alpha$ and $\beta$ such that $|\beta|\le \mu$. 
We can write
$${\rm Op}_h(A_1){\rm Op}_h(A_2)-{\rm Op}_h(A_1A_2)=\sum_{\ell_1=1}^3\sum_{\ell_2=1}^3
\left({\rm Op}_h(A_1^{(\ell_1)}){\rm Op}_h(A_2^{(\ell_2)})-{\rm Op}_h(A_1^{(\ell_1)}A_2^{(\ell_2)})\right)$$
$$=\sum_{\ell_1=1}^2\sum_{\ell_2=1}^2
\left({\rm Op}_h(A_1^{(\ell_1)}){\rm Op}_h(A_2^{(\ell_2)})-{\rm Op}_h(A_1^{(\ell_1)}A_2^{(\ell_2)})\right)$$
 $$+\sum_{\ell_2=1}^2
\left({\rm Op}_h((1-\eta)){\rm Op}_h(A_2^{(\ell_2)})-{\rm Op}_h((1-\eta)A_2^{(\ell_2)})\right)$$
$$+\sum_{\ell_1=1}^2
\left({\rm Op}_h(A_1^{(\ell_1)}){\rm Op}_h((1-\eta)r_0^{1/2})-{\rm Op}_h(A_1^{(\ell_1)}(1-\eta)r_0^{1/2})\right)$$
 $$+\left({\rm Op}_h(1-\eta){\rm Op}_h((1-\eta)r_0^{1/2})-{\rm Op}_h((1-\eta)^2r_0^{1/2})\right).$$
Now we apply Proposition 2.6 to the operators in the double sum in the right-hand side, Proposition 2.5 to the operators
in the next two sums and Proposition 2.2 to the operator in the last term. Thus, taking into account the bounds
 (\ref{eq:6.17}) and (\ref{eq:6.18}), we arrive at (\ref{eq:6.13}) and (\ref{eq:6.16}). The bounds (\ref{eq:6.12}) 
 and (\ref{eq:6.15}) follow from Propositions 2.1 and 2.4.
\eproof

In view of (\ref{eq:6.9}) we have $A_2=-(c_1^2-c_2^2)a_1a_2$. By assumption, $|c_1^2-c_2^2|\ge c>0$, so Proposition 6.2 together
with the estimate (\ref{eq:6.8}) imply
\begin{equation}\label{eq:6.21}
\left\|{\rm Op}_h(A_1A_2)f\right\|_{L^2(\Gamma)}\lesssim h^{1/5}\|f\|_{H_h^1(\Gamma)}+{\cal E}(h,\theta,\tau)\|f\|_{H_h^1(\Gamma)}
\end{equation}
provided $\theta\ge h^{2/5}$, $\theta$ satisfies (\ref{eq:6.2}) and (\ref{eq:6.11}), and $\tau\ge h^{1/2}$, $\tau$ satisfies (\ref{eq:6.14}), where
${\cal E}={\cal E}_2$
if (\ref{eq:1.8}) holds, ${\cal E}={\cal E}_1$ if 
(\ref{eq:1.6}) holds and $z\in Z^+$, and ${\cal E}={\cal E}_3$ if 
(\ref{eq:1.6}) holds and $z\in Z^-$, where 
$${\cal E}_2=h\theta^{-7/2}+h^{\mu/2}\log(h^{-1})\theta^{-d-\mu-1},$$
$${\cal E}_3=h\tau^{-2}+h^{\mu/2}\log(h^{-1})\tau^{-d-\mu-1}.$$
Set
$$\theta_2(h)=
\left\{
\begin{array}{l}
\left(h^{\mu/2}\log(h^{-1})\right)^{1/(d+\mu+1)}
\quad\mbox{if}\quad\mu\le \frac{4}{3}(d+1),\\
h^{2/7}
\quad\mbox{if}\quad \mu>\frac{4}{3}(d+1), 
\end{array}
\right.
$$
and
$$\tau_3(h)=\left(h^{\mu/2}\log(h^{-1})\right)^{1/(d+\mu+1)}.$$
It is easy to check that $\theta\gg\theta_2$ implies ${\cal E}_2\ll 1$ 
and $\theta\ge h^{2/7}$ satisfies (\ref{eq:6.2}) and (\ref{eq:6.11}). Similarly, 
 $\tau\gg\tau_3$ implies ${\cal E}_3\ll 1$ and $\tau\ge h^{1/2}$ satisfies (\ref{eq:6.14}). Since $A_1A_2=(r_0+1)^{1/2}$ and 
$$\left\|{\rm Op}_h\left((r_0+1)^{1/2}\right)f\right\|_{L^2(\Gamma)}\sim\|f\|_{H_h^1(\Gamma)},$$
the estimate (\ref{eq:6.21}) implies 
\begin{equation}\label{eq:6.22}
\|f\|_{H_h^1(\Gamma)}\lesssim h^{1/5}\|f\|_{H_h^1(\Gamma)}+{\cal E}\|f\|_{H_h^1(\Gamma)}.
\end{equation}
Taking $h$ and ${\cal E}$ small enough we deduce from (\ref{eq:6.22}) that $\|f\|_{H_h^1(\Gamma)}=0$, and hence 
$f\equiv 0$, as desired. In other words, the regions $\theta\gg\theta_2$ when (\ref{eq:1.8}) holds, and 
$\theta\gg\theta_1$ and $\tau\gg \tau_3$ when (\ref{eq:1.6}) holds
are eigenvalue-free regions. It is easy to see that these regions correspond
to eigenvalue-free regions in Theorem 1.1 on the $\lambda$ plane.

\end{document}